\documentclass[12pt]{article}
\usepackage[T2A]{fontenc}
\usepackage[cp866]{inputenc}
\usepackage[english,russian]{babel}
\usepackage[tbtags]{amsmath}
\usepackage{amsfonts,amssymb,mathrsfs,amscd,amsthm,hyper}
\usepackage{mathtools}
\RequirePackage{graphicx}
\theoremstyle{plain}

\theoremstyle{plain}
\newtheorem{theorem}{Теорема}
\newtheorem{Proposi}{Предложение}
\newtheorem{propos}{Утверждение-гипотеза}
\newtheorem{lemma}{Лемма}

\theoremstyle{definition}
\newtheorem{definition}{Определение}
\newtheorem{remark}{Замечание}
\theoremstyle{main}
\newtheorem{LemmaN}{Лемма 9}
\newtheorem{LemmaS}{Лемма 7}
\def\({\left(}
\def\){\right)}
\def\Re{\operatorname{Re}}
\def\Im{\operatorname{Im}}
\def\const{\operatorname{const}}
\def\supp{\operatorname{supp}}
\def\mint{\operatorname{int}}
\def\Supp{\operatorname{Supp}}
\def\Rob{\operatorname{Rob}}
\def\NOD{\operatorname{НОД}}
\def\proj{\operatorname{proj}}

\def\LG{\operatorname{LG}}
\def\LG{\operatorname{LG}}

\def\dist{\operatorname{dist}}
\def\sS{S}

\def\fQ{\mathfrak Q}
\def\kk{\mathfrak K}

\def\frame{\operatorname{frame}}
\def\PP{\mathbb P}
\def\mDD{\mathbb D}
\def\sU{\mathscr U}
\def\sH{\mathscr H}
\let\mhat\widehat
\def\NN{\mathbb N}
\def\CC{\mathbb C}
\def\RR{\mathbb R}
\def\ZZ{\mathbb Z}
\def\HH{\mathscr H}

\def\sP{\mathscr P}
\def\sV{\mathscr V}
\def\sN{\mathscr N}
\def\sL{\mathscr L}
\def\KK{\mathscr K}

\def\FF{\mathscr F}
\def\HH{\mathscr H}
\def\GG{\mathscr G}
\def\DD{\mathscr D}
\def\RS{\mathfrak R}
\def\zz{\mathbf z}

\def\balpha{\boldsymbol\alpha}
\def\ba{\mathbf a}

\def\EE{E}
\def\FF{F}
\def\myFF{F}
\def\mdeg{\operatorname{deg}}
\def\mcap{\operatorname{cap}}
\let\myo\overline\let\myh\widehat\let\eps\varepsilon\let\pfi\varphi
\let\le\leqslant\let\leq\leqslant\let\geq\geqslant
\let\myt\widetilde

\let\le=\leqslant

\def\bad{\spaceskip=0.33emplus0.6emminus0.15em\immediate\write5{\string\bad}}
\begin{document}

\title{On the distribution of zeros of the Hermite--Pade polynomials for three
algebraic functions $1,f,f^2$ and the global topology of the Stokes lines for
some differential equations of the third order}

\author{Sergey Suetin}


\date{26.12.2013}

\maketitle

\begin{abstract}
The paper presents some heuristic results about the distribution of zeros
of Hermite--Pade polynomials of first kind for the case of three
functions $1,f,f^2$, where $f$ has the form
$f(z): = \prod\limits_ {j = 1 } ^3 (z-a_j) ^ {\alpha_j} $,
$\alpha_j \in \mathbb C\setminus \mathbb Z $,
$ \sum \limits_ {j = 1 } ^ 3 \alpha_j = 0 $, $ f (\infty) = 1 $ .
Answers are given in terms related to the problem of extreme
minimum capacity of a plane Nuttall condenser, two plates of which intersect
(``hooked '' for each other) in a five 5
points: the branch points $ a_1, a_2, a_3 $ of function $ f $, at $ v_1 = v $,
where $ v = v (a_1, a_2, a_3) $ the Chebotarev point corresponding to triple
points $ a_1, a_2, a_3 $, and another ``unknown'' at $ v_2 \neq v_1, a_1, a_2, a_3 $ ( see
Fig1--Fig3).

The connection between the distribution of zeros of Hermite--Pade polynomials
and global topology of the Stokes lines and the asymptotic behavior of
Liouville--Green solutions of a class of homogeneous linear differential
equations of the third order containing a large parameter is the free term
is established.

The basic idea of the new and still heuristic approach is to reduce at first
some theoretical potential vector equilibrium problem to the scalar problem
with the external field, and then use the general method  of Gonchar--Rakhmanov
developed in 1987 for solving the Varga problem `` about $ 1/9 $''.

It is supposed that in the general case of an arbitrary algebraic function $ f
$ it is imposible to constract the Nuttall condenser for a set of three
functions $ 1, f, f ^ 2 $ without the knowledge of the structure of the Stahl
compact for function $ f $. Namely, the Nuttall condenser is constructed using
Green's function $ g_ {D} (z, \infty) $ for the Stahl domain $ D $ and ``
core'' of Stahl compact, which consists of `` effective'' branch points of $ f
$ and corresponding Chebotarev points.

Библиография: 71 название.
\end{abstract}

\markright{Об асимптотике полиномов Эрмита--Паде}

\section{Введение и постановка задачи}\label{s1}

\subsection{}\label{s1s1}
В настоящей работе установлена связь между распределением нулей полиномов
Эрмита--Паде для системы из трех функций $1,f,f^2$ (см.~\eqref{1.5}), где
\begin{equation}
f(z):=\prod_{j=1}^3(z-a_j)^{\alpha_j},\qquad 2\alpha_j\in\CC\setminus\ZZ,
\quad \sum_{j=1}^3\alpha_j=0,\quad f(\infty)=1,
\label{01.1}
\end{equation}
$a_1,a_2,a_3$ -- три различных точки комплексной плоскости $\CC$, находящиеся
в ``общем положении'', и глобальной структурой линий Стокса и $\LG$-асимптотикой
для специального линейного однородного дифференциального уравнения 3-го порядка,
содержащего большой параметр при свободном члене (см.~\eqref{ma1.8}).

Вопрос о распределении нулей таких полиномов исследуется с помощью
общего метода Гончара--Рахманова, развитого в~\cite{GoRa87} при решении задачи
Варги~``об $1/9$'' (см.~\cite{Var74}). Ответ дается в терминах, связанных с решением некоторой
экстремальной задачи о минимуме (обобщенной) емкости плоского
конденсатора, состоящего из 2-х пластин, соприкасающихся в точках $a_j$
и находящихся под действием внешнего поля. Две пластины {\it экстремального}
конденсатора~$(E,F)$ оказываются
зацепленными в пяти точках: точках ветвления $a_1,a_2,a_3$ функции $f$ и еще двух точках
$v_1,v_2$ -- некоторых трансцендентных параметрах задачи; при этом точка $v_1$
совпадает с классической точкой Чеботарёва $v$ для трех точек $a_1,a_2,a_3$, а $v_2\neq
v_1,a_1,a_2,a_3$ (см. рис.~\ref{Fig1}--\ref{Fig3}). Пластины
соответствующего экстремального конденсатора обладают определенным свойством
``симметрии'' (или $\sS$-свойством), вполне аналогичным введенному ранее
в~\cite{GoRa87}, и состоят из замыканий критических траекторий квадратичного
дифференциала $\fQ(z)\,dz^2$.
В настоящей работе
устанавливается, что предельное
распределение нулей полиномов Эрмита--Паде существует и совпадает с равновесной
мерой $\lambda_F$ для пластины $F$ экстремального конденсатора $(E,F)$. Аналогично
с помощью этого же метода устанавливаются формулы слабой асимптотики для
функции остатка и отношения полиномов Эрмита--Паде. Все эти предельные
соотношения носят характер ``сходимости по емкости''. Однако при этом
оказывается, что такая сходимость обладает замечательным свойством: {\it
она выдерживает дифференцирование любое конечное число раз}. Затем, с целью
вывода формул сильной асимптотики, доказывается, что функция остатка и полиномы
Эрмита--Паде удовлетворяют
линейному однородному дифференциальному уравнению третьего порядка
с полиномиальными коэффициентами фиксированной степени и большим параметром
при свободном члене. Устанавливается, что соответствующее характеристическое
уравнение 3-й степени обладает тем свойством, что его {\it трехлистная риманова
поверхность обладает каноническим разбиением на листы
в смысле Наттолла} (см.~\cite{KovSu14} и ниже
п.~\ref{s6s1}). Поэтому линии Стокса для полученного
дифференциального уравнения -- это в точности аналитические дуги, составляющие
пластины экстремального конденсатора $E\cup F$, а поступательные пути
(``progressive path'') -- сопряженные траектории квадратичного дифференциала
$\fQ(z)\,dz^2$, которые следует рассматривать на {\it трех различных листах
характеристической римановой поверхности} и тогда они не будут пересекаются.
После этого оказывается возможным предъявить в явном виде три асимптотически
независимых $\LG$-приближения $\myt{w}_{n,1},\myt{w}_{n,2},\myt{w}_{n,3}$,
подходящая линейная комбинация которых позволяет асимптотически приблизить
любое наперед заданное решение дифференциального уравнения. Таким образом, для вывода формул
сильной асимптотики для полиномов Эрмита--Паде и функции остатка используется
обобщение классического метода Лиувилля--Грина\footnote{Как отмечено в
монографии Сегё~\cite[глава 6]{Sze62} в асимптотической теории ортогональных
многочленов этот метод принято называть методом Лиувилля--Стеклова; см.
также~\cite{Ste07},~\cite{Nev84}.},
основанного на асимптотическом приближении
произвольного решения дифференциального уравнения 2-го порядка некоторой
линейной комбинацией 2-х $\LG$-приближений. Для вывода соответствующих формул
типа Биркгофа выясняется глобальная структура линий Стокса и соответствующих
канонических путей. Делается это следующим образом. После того как установлено
существование экстремального конденсатора и показано, что обе его пластины
$E$ и $F$ состоят из замыканий критических траекторий квадратичного
дифференциала $\fQ(z)\,dz^2$ и обладают $\sS$-свойством, методом, предложенным
А.~А.~Гончаром в 1990-х годах в статьях~\cite{GoRaSu91},~\cite{GoRaSu92}, строится
трехлистная риманова поверхность~$\RS_3$ (см. также~\cite{KovSu14}).
Оказывается, что на этой поверхности $\RS_3$ существует абелев интеграл
$\pfi(\zz)$, $\zz\in\RS_3$, третьего рода с чисто мнимыми периодами, такой, что
$\Re\pfi(\zz)$ -- однозначная гармоническая, $\pfi'(\zz)$ -- однозначная
мероморфная функции на $\RS_3$,
\begin{equation}
\fQ(z)\,dz^2=-(\pfi'(z^{(2)}))^2\,dz^2
\label{01.2}
\end{equation}
(тем самым квадратичный дифференциал поднимается с $\myo\CC$ на $\RS_3$) и с помощью $\pfi$
можно естественным образом определить глобальное разбиение $\RS_3$ на ``листы''
$z^{(j)}$, $j=1,2,3$:
\begin{equation}
\Re\pfi(z^{(3)})>\Re\pfi(z^{(2)})>\Re\pfi(z^{(1)}),\quad z\notin E\cup F.
\label{01.3}
\end{equation}
Линии Стокса для дифференциального уравнения~\eqref{ma1.8} надо теперь
рассматривать как критические траектории
\begin{equation}
-(\pfi'(\zz))^2\,dz^2>0
\label{01.4}
\end{equation}
на римановой поверхности, соответствующие линии уровня $\Re\pfi(\zz)=\const$
задают ``вертикальные'' линии, непересекающиеся на $\RS_3$ (но пересекающиеся
после их проектирования на $\myo\CC$). Непересекающиеся на $\RS_3$ ``горизонтальные''
линии $\Im\pfi(\zz)=\const$, соответствующие сопряженным траекториям
$$
-(\pfi'(\zz))^2\,dz^2<0,
$$
определяют канонические пути, необходимые для построения $\LG$-теории. Через
каждую точку $\zz$ рп $\RS_3$ проходит в точности одна вертикальная линия и
одна горизонтальная линия. Разумеется, после проектирования точки $\zz\in\RS_3$
на $\myo\CC$ эта единственность теряется.
Благодаря \eqref{01.3} три функции
\begin{equation}
\myt{w}_{n,j}(z)
=\exp\biggl\{n\pfi(z^{(j)})-\int_{a_1}^{z^{(j)}}\sum_{k=1,k\neq{j}}^3
\frac{\pfi''(z^{(j)})}{\pfi'(z^{(j)})-\pfi'(z^{(k)})}\,dz\biggr\},
\quad
j=1,2,3,
\label{01.5}
\end{equation}
имеют разную асимптотику в открытом множестве $\myo{\CC}\setminus(E\cup F)$ и их можно
взять в качестве $\LG$-приближений. Любое решение дифференциального
уравнения~\eqref{ma1.8} оказывается возможным приблизить линейной комбинацией вида
\begin{equation}
C_1\myt{w}_{n,1}(z)+C_2\myt{w}_{n,2}(z)+C_3\myt{w}_{n,3}(z),
\label{01.6}
\end{equation}
в которой в силу неравенств~\eqref{01.3} всегда есть только одно доминирующее
при $n\to\infty$ слагаемое. Выбор этого доминирующего слагаемого зависит
от того, какому именно листу рп $\RS_3$ соответствуют заданные начальные
условия для дифференциального уравнения~\eqref{ma1.8}.

Таким образом, в настоящей работе установлено существование аналитических
кривых, на которых в пределе (при $n\to\infty$) распределяются нули полиномов
Эрмита--Паде. Оказывается, что эти кривые $F$ состоят из замыканий критических
траекторий квадратичного дифференциала. Если теперь включить в рассмотрение
и вторую часть $E$ критических траекторий этого квадратичного дифференциала и
с помощью этой пары $(E,F)$ построить трехлистную риманову поверхность, то на
этой римановой поверхности естественным образом возникает разбиение на листы
в смысле Наттолла. Сопряженные траектории для
квадратичного дифференциала окажутся поступательными путями, необходимыми
для построения аналога классической $\LG$-теории для дифференциального
уравнения третьего порядка.

Функция $f$ класса~\eqref{01.1} зависит от точек $a_1,a_2,a_3$ как от некоторых
параметров. От этих же параметров будет зависеть и соответствующее
дифференциальное уравнение 3-го порядка. Тем самым предложенный здесь подход
построения трехмерного аналога классической $\LG$-теории оказывается пригодным
для целого {\it семейства} дифференциальных уравнений 3-го порядка с большим
параметром при свободном члене. Этот класс можно расширить, если
вместо~\eqref{01.1} рассмотреть функции вида
$$
f(z):=\prod_{j=1}^p(z-a_j)^{\alpha_j},\qquad\alpha_j\in\CC\setminus\ZZ,
\quad \sum_{j=1}^p\alpha_j=0,\quad f(\infty)=1,
$$
где $p\geq3$, $a_1,a_2,\dots,a_p$ -- попарно различные точки комплексной
плоскости $\CC$, находящиеся в ``общем положении''
(см.~\cite{MaRaSu12},~\cite{MaRaSu13}).

Подчеркнем, что основные результаты работы, относящиеся к асимптотическим
формулам для полиномов Эрмита--Паде, носят эвристический характер и
их следует рассматривать как гипотезы, нуждающиеся в строгом обосновании.
Фактически здесь предлагается некоторая {\it программа} по исследованию на
основе метода Гончара--Рахманова
(см.~\cite{GoRa87},~\cite{Gon03},~\cite{AptLopSaf11})
структуры линий Стокса, а затем на основе метода Лиувилля--Стеклова
построению глобальной
асимптотической
теории для некоторого класса линейных однородных дифференциальных уравнений
третьего порядка с полиномиальными коэффициентами и большим параметром при
свободном члене.

Отметим следующий важный факт, вытекающий из вышеприведенных результатов.
Для классического случая дифференциального уравнения 2-го порядка структура
линий Стокса описывается в терминах {\it квадратичного} дифференциала, а
характеристическое уравнение для этого случая -- алгебраическое уравнение
второго порядка, т.е. также {\it квадратичное}. Таким образом имеется некоторая
априорная согласованность между ``управляющим'' дифференциалом и характеристическим
уравнением: оба они являются квадратичными.
Вполне естественно, что в этом случае глобальную асимптотическую
$\LG$-теорию
оказывается возможным построить практически для {\it любого} дифференциального
уравнения 2-го порядка (см.~\cite{Fed65}, \cite{EvgFed66},~\cite{Fed86},
\cite{EsiSha10},~\cite{GalSha06}), даже и нелинейного~\cite{Fed86b},~\cite{Ste02}.

Совсем иначе обстоит дело в случае дифференциального уравнения третьего
порядка (см.~\cite{KaKoTa10},~\cite{KaKoTa11}).
За структуру линий Стокса и в этом случае отвечает {\it квадратичный}
дифференциал, но характеристическое уравнение -- алгебраическое уравнение 3-й
степени, т.е. {\it кубическое}. В этом случае глобальную асимптотическую
$\LG$-теорию оказывается возможным построить только в таком классе
дифференциальных уравнений, когда этот квадратичный дифференциал и кубическое
уравнение оказываются {\it согласованными} друг с другом. Это оказывается
справедливо для случая, когда соответствующая характеристическая трехлистная
риманова поверхность допускает каноническое разбиение на листы
(см.~\eqref{01.3}, а также п.~\ref{s6s1}).

Результаты настоящей работы делятся на две части. Первая из них (в основном это
результаты вспомогательного и подготовительного характера; см.
леммы~\ref{lem2}--\ref{lemma2})
состоит из полностью и строго доказанных результатов. Вторая основная часть
результатов настоящей работы состоит из утверждений (см.
утверждения-гипотезы~\ref{pro1}--\ref{pro4}),
полученных с помощью нестрогих эвристических рассуждений, эти утверждения
следует рассматривать как гипотезы, которые нуждаются в последующем строгом
обосновании.

Отметим, что функция $f$ удовлетворяет дифференциальному уравнению первого
порядка
$$
\frac{f'}{f}(z)=\sum_{j=1}^3\frac{\alpha_j}{z-a_j}=\frac{B(z)}{A(z)},
$$
$B\in\PP_1$.

\subsection{}\label{s1s2}
Для того чтобы прояснить наиболее существенные детали нашего эвристического
подхода, проиллюстрируем его сначала на примере полиномов Эрмита--Паде для пары
функций $1,f$, где как и прежде $f$ задана представлением~\eqref{01.1}. А именно,
пусть нетривиальные полиномы $P_{n,0},P_{n,1}\not\equiv0$ степени $\leq{n}$
удовлетворяют соотношению
\begin{equation}
T_n(z):=(P_{n,0}+P_{n,1}f)(z)=O\(\frac1{z^{n+1}}\),\quad z\to\infty.
\label{01.7}
\end{equation}
Ясно, что $P_{n,0},P_{n,1}$ -- полиномы Паде для функции $f$, $T_n$ --
соответствующая функция остатка, $[n/n]_f:=-P_{n,0}/P_{n,1}$ -- диагональная
аппроксимация Паде функции $f$ порядка $n$.

Вопрос о распределении нулей
полиномов Паде и асимптотике отношения $P_{n,0}/P_{n,1}$ решен Шталем для
достаточно широкого класса функций, которому принадлежат и функции
вида~\eqref{01.1}. Изложим кратко основные положения этой теории
Шталя~\cite{Sta97b} в применении к функциям вида~\eqref{01.1}.

Пусть $\kk_f$ -- семейство всех компактов $K\not\ni\infty$ не разбивающих
комплексную плоскость и таких, что функция $f$ продолжается из окрестности
бесконечно удаленной точки $z=\infty$ в дополнение к $K$ как голоморфная
(однозначная аналитическая) функция, $f\in\sH(\myo\CC\setminus{K})$. Компакты
семейства $\kk_f$ будем называть допустимыми (для заданной функции~$f$). Имеет
место следующая
\begin{theorem}[(Г.~Шталь,~\cite{Sta97b})]\label{te1}
Существует единственный допустимый компакт $S$, $S\in \kk_f$, такой, что
\begin{equation}
\mcap{S}=\inf_{K\in \kk_f}\mcap{K}.
\label{01.8}
\end{equation}
Компакт $S$ не разбивает плоскость, состоит из трех аналитических
дуг $S_1,S_2,S_3$ -- замыканий критических траекторий квадратичного дифференциала
\begin{equation}
\fQ(z)\,dz^2=-\frac{z-v}{A_3(z)}\,dz^2,\quad A_3(z)=\prod_{j=1}^3(z-a_j),
\label{01.9}
\end{equation}
и обладает следующим свойством симметрии ($\sS$-свойством):
\begin{equation}
\frac{\partial g_S(\zeta,\infty)}{\partial n^{+}}=
\frac{\partial g_S(\zeta,\infty)}{\partial n^{-}},\qquad \zeta\in {S^0},
\label{01.10}
\end{equation}
где $g_S(\zeta,\infty)$ -- функция Грина для области Шталя
$D:=\myo\CC\setminus{S}$, ${S^0}$ -- объединение открытых дуг $S^0_1\cup
S^0_2\cup S^0_3$, замыкания которых составляют $S$, $\partial/\partial n^{\pm}$ -- производные
по внутренней по отношению к области $D$ нормали в точке $\zeta$, взятые с
разных сторон~$S$.
\end{theorem}

\begin{remark}[(см.~\cite{Kuz80})]\label{rem1}
Компакт $S$ зависит только от взаимного расположения точек $a_1,a_2,a_3$,
точка Чеботарёва $v$ лежит внутри треугольника с вершинами в точках
$a_1,a_2,a_3$. В дальнейшем мы всегда считаем, что точки $a_1,a_2,a_3$
находятся в ``общем положении''. Это в частности означает, что индекс $n$
является нормальным для полиномов Паде: $\mdeg{P_{n,j}}=n$, $j=0,1$, и для
функции остатка имеем:
$$
T_n(z)=\frac{C_n}{z^{n+1}}+O\(\frac1{z^{n+2}}\),\quad C_n\neq0.
$$
Те же самые соглашения относятся в дальнейшем и к полиномам Эрмита--Паде.
\end{remark}

\subsection{}\label{s1s2s2}
Пусть $\lambda_S$ -- единичная равновесная мера для компакта $S$:
$V^{\lambda_S}(z)\equiv\gamma_S$, $z\in{S}$,
$V^{\lambda_S}(z)=-\int\log|z-\zeta|\,d{\lambda_S}(\zeta)$ -- логарифмический
потенциал, $\gamma_S$ -- постоянная Робена. Поскольку
$g_S(\zeta,\infty)=\gamma_S-V^{\lambda_S}(z)$, соотношение~\eqref{01.10}
эквивалентно равенству
\begin{equation}
\frac{\partial V^{\lambda_S}(\zeta)}{\partial n^{+}}=
\frac{\partial V^{\lambda_S}(\zeta)}{\partial n^{-}},\qquad \zeta\in{S^0}.
\label{01.11}
\end{equation}
Справедлива следующая

\begin{theorem}[(Г.~Шталь,~\cite{Sta97b})]\label{te2}
Пусть $P^{*}_{n,j}(z)=z^n+\dots$ -- монический полином, соответствующий
полиному $P_{n,j}$, $j=0,1$, $T^{*}_n(z)=1/z^{n+1}+\dots$ -- нормированная
функция остатка. Тогда
\begin{gather}
|P^{*}_{n,j}(z)|^{1/n}\overset\mcap\longrightarrow e^{-V^{\lambda_S}(z)},
\quad n\to\infty,\quad z\in D,\quad j=0,1,
\label{01.12}\\
\frac{P_{n,0}}{P_{n,1}}(z)\overset\mcap\longrightarrow -f(z),
\quad n\to\infty,\quad z\in D,
\label{01.13}\\
|T^{*}_{n}(z)|^{1/n}\overset\mcap\longrightarrow e^{V^{\lambda_S}(z)},
\quad n\to\infty,\quad z\in D.
\label{01.14}
\end{gather}
\end{theorem}
Соотношения~\eqref{01.12}--\eqref{01.14} понимаются в смысле сходимости по
емкости; \eqref{01.12} и~\eqref{01.14} -- это асимптотика корня $n$-й степени,
\eqref{01.13} -- асимптотика отношения; вопрос о сильной асимптотике обсуждается
ниже в~\S~\ref{s5}.

Для произвольного полинома $Q\in\CC[z]$ положим
$\mu(Q)=\sum_{\zeta:Q(z)=0}\delta_\zeta$ -- мера, считающая нули этого
полинома. Тогда соотношения~\eqref{01.12} эквивалентны соотношениям
\begin{equation}
\frac1n\mu(P_{n,j})\to\lambda_S,\quad n\to\infty,\quad j=0,1;
\label{01.15}
\end{equation}
соотношения~\eqref{01.15} понимаются в смысле слабой сходимости в пространстве
мер.

Отметим очень важное свойство сходимости по емкости:
{\it соотношения~\eqref{01.12}--\eqref{01.14} можно дифференцировать любое конечное
число раз}. А именно, имеет место следующая лемма.

\begin{lemma}\label{lem2}
В условиях теоремы Шталя справедливы следующие соотношения
\begin{gather}
\frac1n\frac{P'_{n,j}}{P_{n,j}}(z)\overset\mcap\longrightarrow
\myh\lambda_S(z),\quad
\frac1{n^2}\frac{P''_{n,j}}{P_{n,j}}(z)\overset\mcap\longrightarrow
\myh\lambda_S(z)^2, \quad n\to\infty,\quad z\in D,\quad j=0,1,
\label{01.16}\\
\frac1n\frac{R'_{n}}{R_{n}}(z)\overset\mcap\longrightarrow
-\myh\lambda_S(z),\quad
\frac1{n^2}\frac{R''_{n}}{R_{n}}(z)\overset\mcap\longrightarrow
\myh\lambda_S(z)^2, \quad n\to\infty,\quad z\in D.
\label{01.17}
\end{gather}
\end{lemma}
В~\eqref{01.16}--\eqref{01.17}
$$
\myh\lambda_S(z)=\int\frac{d\lambda_S(\zeta)}{z-\zeta}
$$
-- преобразование Коши равновесной меры $\lambda_S$,
$\myh\lambda_S(z)=-\sV'_{\lambda_S}(z)$, где
$$
\sV_{\lambda_S}(z)=\int\log\frac1{z-\zeta}\,d\lambda_S(\zeta)
$$
-- комплексный потенциал меры $\lambda_S$.

\subsection{}\label{s1s0}
Отметим следующую связь результата Шталя (см. теорему~\ref{te1}) о компакте
минимальной емкости с теоретико-потенциальной задачей о {\it о максимуме
равновесной энергии} (см.~\cite{PeRa94}).

Пусть $K\in\kk_f$ -- допустимый компакт для функции~$f$, $M_1(K)$ --
пространство всех единичных (положительных борелевских) мер с носителями
на~$K$, $V^\mu(z)=-\int\log|z-\zeta|\,d\mu(\zeta)$ -- логарифмический потенциал
меры $\mu\in M_1(K)$. В классе $M_1(K)$ существует единственная (в этом классе)
{\it равновесная} мера $\lambda_K\in M_1(K)$: $V^{\lambda_K}(z)\equiv\gamma_K$,
$z\in K$, $\gamma_K$ -- постоянная Робена для $K$, имеем
$e^{-\gamma_K}=\mcap{K}$.

Мера $\lambda_K$ и только эта мера является {\it экстремальной} мерой для
функционала энергии $J(K;\mu)$,
\begin{equation}
J(K;\mu):=\int_K V^\mu(z)\,dz=\iint_{K\times K}\log\frac1{|z-\zeta|}\,d\mu(z)\,d\mu(\zeta),
\label{t1}
\end{equation}
в классе $M_1(K)$:
\begin{equation}
J(K;\lambda_K)=\min_{\mu\in M_1(K)}J(K;\mu).
\label{t2}
\end{equation}
При этом $J(K;\lambda_K)=\gamma_K$. Поскольку
$e^{-\gamma_K}=\mcap{K}$, то задача
\begin{equation}
\min_{K\in\kk_f}\mcap{K}
\label{t3}
\end{equation}
эквивалентна ``$\max$-$\min$''-задаче:
\begin{equation}
\max_{K\in\kk_f}J(K;\lambda_K)
=\max_{K\in\kk_f}\min_{\mu\in M_1(K)}J(K;\mu).
\label{t4}
\end{equation}
Из вышесказанного вытекает, что если $F\in\kk_f$ -- допустимый компакт, для
которого
$$
J(F;\lambda_F)=\max_{K\in\kk_f}J(K;\lambda_K),
$$
то $F=S\in\sS$, т.е. (единственный) экстремальный компакт необходимо обладает
$\sS$-свойством (ср. утверждение~\ref{pro1}). При этом $J(F;\lambda_F)=-\log{\mcap{S}}$.


\section{Формулировка основных результатов}\label{s2}

\subsection{}\label{s2s1}
Дадим определение полиномов Эрмита--Паде первого рода для трех функций
$1,f,f^2$, где функция $f$ имеет вид~\eqref{01.1}.

Для системы из трех функций $1,f,f^2$
и произвольного $n\in\NN$ полиномы Эрмита--Паде I~рода
(см.~\cite{GoRa81}, \cite{Nut81}, \cite{Nut84}, \cite{NuTr87},
\cite{ApLy10},~\cite{ApLyTu11},~\cite{Rak11},~\cite{Sor10})
$Q_{n,0},Q_{n,1},Q_{n,2}\not\equiv0$ степеней $\leq{n}$,
определяются (неоднозначно) из соотношения
(ср.~\eqref{01.7}):
\begin{equation}
\myt{R}_n(z):=
(Q_{n,0}\cdot1+Q_{n,1}f+Q_{n,2}f^2)(z)=O\(\frac1{z^{2n+2}}\),\qquad z\to\infty;
\label{1.5}
\end{equation}
естественно, для невырождения\footnote{Легко
видеть, что в силу условия $\sum\limits_{j=1}^3\alpha_j=0$
при $\alpha_1=\pm1/2$ имеем $\alpha_{2,3}\neq\pm1/2$.
Поэтому в заданном классе функций~\eqref{01.1} полного вырождения никогда не
наступает, ответ все равно существует,
но при $\alpha_1=\pm1/2$ он будет уже {\it другой\/}.}
задачи теперь нужно дополнительное условие на $f$: все $\alpha_j\neq\pm1/2$.

Пусть $Q^{*}_{n,j}$, $j=0,1,2$, -- соответствующие монические полиномы.
Рассмотрим следующую задачу: {\it найти распределение нулей полиномов
$Q_{n,j}$}, т.е. найти при $n\to\infty$ пределы нормированных считающих мер
$$
\frac1n\mu(Q_{n,0}),\quad
\frac1n\mu(Q_{n,1}),\quad\frac1n\mu(Q_{n,2}).
\label{1.6}
$$
Ясно, что эта задача эквивалентна задаче об асимптотике корня $n$-й степени
для монических полиномов:
$$
|Q^{*}_{n,j}|^{1/n},\quad n\to\infty,\quad z\in\CC,\quad j=0,1,2.
$$
Введем в рассмотрение семейство {\it допустимых} компактов $\GG=\{\Gamma\}$ и
{сопряженное} ему семейство компактов $\GG^{\perp}=\{\Gamma^{\perp}\}$. Компакт
$\Gamma\in\GG$, если $\Gamma=\bigcup_{j=1}^3\Gamma_j$, $\Gamma_j$ --
аналитическая дуга, соединяющая точки $a_j$ и $v$ и гомотопная дуге $S_j$
($\Gamma_j$ пересекаются
только по точке $v$). Дуальный компакт $\Gamma^{\perp}$ не пересекается с $\Gamma$
и состоит из $6$ дуг,
пересекающихся в (``исходящих'' из) произвольной точке $c\notin{S}$. Три из них
заканчиваются в точке Чеботарёва $v$, три других -- в точках ветвления $a_j$
(см. рис.~\ref{Fig1}; это описание неформальное, можно дать и более строгое
определение в терминах ``мембран'' и ``хвостов''~\cite{ApTu12b}, см. также
ниже замечание~\ref{rem1.1}). Таким образом, $\myo{\CC}\setminus{\Gamma^{\perp}}$ состоит из {\it трех} непересекающихся
областей $D_j\supset S^0_j$, таких, что $D_j\cap{\Gamma}=\Gamma^0_j$ (фактически области
$D_j$ ``нанизаны'' на дуги $S^0_j$, см. рис. ~\ref{Fig3}).

При фиксированном
компакте $\Gamma$ дуальных ему компактов $\Gamma^{\perp}$ много. Пусть $M_1(\Gamma)$
-- единичные меры с носителями на $\Gamma$, $M_1(\Gamma^{\perp})$ -- единичные меры с
носителями на $\Gamma^{\perp}$, $A=\begin{pmatrix}1&-1\\-1&4\end{pmatrix}$ -- матрица
взаимодействия, $\psi(z)=3g_{S}(z,\infty)$ -- внешнее поле на $\Gamma^{\perp}$,
$g_{S}(z,\infty)$ -- функция Грина для области Чеботарёва--Шталя $D_S$. Для мер
$\mu_1\in M_1(\Gamma)$, $\mu_2\in M_1(\Gamma^{\perp})$ определим взаимную $A$-энергию
в поле $\psi(z)$ (поле $\psi$ действует только на $\Gamma^{\perp}$):
\begin{align}
J_{A,\psi}(\mu_1,\mu_2):&=(A\vec{\mu},\vec{\mu})
+2\int_{\Gamma^{\perp}}\psi(\zeta)\,d\mu_2(\zeta)\notag\\
&=\sum_{i,j=1}^2 a_{ij}[\mu_i,\mu_j]+2\int_{\Gamma^{\perp}}\psi(\zeta)\,d\mu_2(\zeta),
\label{2.1}
\end{align}
где $\vec{\mu}=(\mu_1,\mu_2)$,
$$
[\mu,\nu]=\iint\log\frac1{|z-\zeta|}d\mu(z)\,d\nu(\zeta)
$$
-- стандартная (удвоенная) энергия взаимодействия двух мер
(см.~\cite{Lan66},~\cite{GoRa81}, \cite{GoRa85},~\cite{SaTo97}).

При фиксированных $\Gamma\in\GG_f=\GG$ и
$\Gamma^{\perp}\in\GG^{\perp}_f=\GG^{\perp}$
рассмотрим следующую задачу:
$$
\inf_{\substack{\mu_1\in M_1(\Gamma),\\\mu_2\in M_1(\Gamma^{\perp})}}
J_{A,\psi}(\mu_1,\mu_2)=m=m(\Gamma,\Gamma^{\perp}).
$$

Стандартным образом (см.~\cite{Lan66},~\cite{GoRa81},~\cite{SaTo97})
доказывается следующая

\begin{lemma}\label{le0}
Существует (единственная) пара мер $(\lambda_1,\lambda_2)$ такая,
что $\supp\lambda_1=\Gamma\in\GG_f$, $\supp\lambda_2=\Gamma^{\perp}\in\GG^{\perp}_f$ и
$$
J_{A,\psi}(\lambda_1,\lambda_2)=m.
$$
\end{lemma}

\subsection{}\label{s2s2}

Рассмотрим теперь следующую {\it экстремальную задачу}:
\begin{equation}
\sup_{\Gamma\in\GG_f,\Gamma^{\perp}\in\GG^{\perp}_f}m(\Gamma,\Gamma^{\perp})=M.
\label{2.2}
\end{equation}

\begin{propos}\label{pro1}
Пусть функция $f$ -- вида~\eqref{01.1}, где все
$\alpha_j\neq\pm1/2$. Тогда

1) Существует единственная пара (стационарных для функционала
энергии~\eqref{2.1}) компактов $(\EE,\myFF)$, $\EE\in\GG_f$, $\myFF\in\GG^{\perp}_f$,
такая, что
\begin{equation}
m(\EE,\myFF)=M.
\label{ext1}
\end{equation}

2) Для соответствующих равновесных мер $\lambda_E,\lambda_F$,
$\supp\lambda_E=\EE$, $\supp\lambda_F=\myFF$,
имеем:
\begin{equation}
\begin{aligned}
V^{\lambda_E}(z)-V^{\lambda_F}(z)&\equiv\gamma_E,\quad z\in \EE,\\
-V^{\lambda_E}(z)+4V^{\lambda_F}(z)+3g_S(z,\infty)&\equiv \gamma_F,\quad z\in \myFF,
\end{aligned}
\label{2.3}
\end{equation}
при этом $\supp{\lambda_E}=\EE$, $\supp{\lambda_F}=\myFF$ {\rm(т.е. явления
сталкивания нет)}, компакты $E$ и $F$ состоят из конечного числа аналитических
дуг (траекторий квадратичного дифференциала) и пара $(\EE,\myFF)$
обладает $\sS$-свойством:
\begin{align}
\frac{\partial(V^{\lambda_E}-V^{\lambda_F})}{\partial n_{+}}(\zeta)
&=\frac{\partial(V^{\lambda_E}-V^{\lambda_F})}{\partial n_{-}}(\zeta),
\quad\zeta\in\EE^0,
\label{2.4.1}\\
\frac{\partial(4V^{\lambda_F}-V^{\lambda_E}+3g_S(\cdot,\infty))}
{\partial n_{+}}(\zeta)
&=\frac{\partial(4V^{\lambda_F}-V^{\lambda_E}+3g_S(\cdot,\infty))}
{\partial n_{-}}(\zeta),\quad\zeta\in\FF^0.
\label{2.4.2}
\end{align}
\end{propos}

\begin{definition}\label{def1}
Пару $(E,F)=(E;F,\psi)$, удовлетворяющую условиям~\eqref{2.3}--\eqref{2.4.2},
назовем {\it ($\psi$-оснащенным) конденсатором Наттолла}.
\end{definition}

\begin{propos}\label{pro2}
Для полиномов Эрмита--Паде $Q_{n,j}$ (см.~\eqref{1.5}) имеем:
\begin{equation}
\lim_{n\to\infty}\frac1n\mu(Q_{n,0})
=\lim_{n\to\infty}\frac1n\mu(Q_{n,1})
=\lim_{n\to\infty}\frac1n\mu(Q_{n,2})
=\lambda_F.
\label{2.5}
\end{equation}
\end{propos}

\begin{propos}\label{pro3}
Пусть $D^{*}_j$, $j=1,2,3$, -- экстремальные области для задачи~\ref{ext1},
составляющие открытое множество $D^*=D^{*}_1\sqcup D^{*}_2\sqcup D^{*}_3$, где
$D^{*}:=\myo\CC\setminus{F}$ (область $D^{*}_j$ ``нанизана'' на
соответствующую  дугу Шталя $S^{0}_j$; см. рис.~\ref{Fig3}). Положим
$f_j|_{S^{0}_j}=(f^{+}+f^{-})(\zeta)$, $\zeta\in S^{0}_j$, $j=1,2,3$,
$f_j\in\HH(D^{*}_j)$ -- голоморфное продолжение этой функции с дуги $S^{0}_j$
в область $D^{*}_j$. Тогда
\begin{equation}
\frac{Q_{n,1}}{Q_{n,2}}(z)\overset\mcap\longrightarrow-f_j(z),
\quad n\to\infty,\quad z\in D^{*}_j,\quad j=1,2,3.
\label{2.5.2}
\end{equation}
\end{propos}

\begin{remark}\label{rem1.1}
Семейство $\GG^{\perp}$ можно сделать исходным (и тогда определения станут
более строгими). Начнем с открытого множества $D\in\DD_f$: $D=D_1\sqcup
D_2\sqcup D_3$, где $D_j$ -- {\it неналегающие} области такие, что
$D_j\supset{S^0_j}$, $j=1,2,3$ (тем самым, $\partial D_j\ni a_j,v$, каждая
область $D^{*}_j$ ``нанизана'' на соответствующую дугу Шталя $S^{0}_j$).
Тогда семейство $\{\partial D:D\in\DD_f\}$ определяет семейство $\GG^{\perp}_f$. Поскольку
существует дуга $\Gamma_j\subset D_j$, соединяющая точки $v$ и $a_j$, то компакты
$\Gamma=\Gamma_1\sqcup\Gamma_2\sqcup\Gamma_3$ задают семейство $\GG_f$.
Вместо пары соотношений равновесия~\eqref{2.3} имеем одно (эквивалентное им)
соотношение равновесия
\begin{equation}
3V^{\lambda}(z)+G^{\lambda}_{\EE}(z)+3g_S(z,\infty)\equiv\gamma_F,
\quad z\in \myFF,
\label{2.6}
\end{equation}
$G_{\EE}^\lambda$ -- гринов (относительно $\EE$) потенциал меры $\lambda\in
M_1(F)$ с носителем
на $\myFF$, компакт $\myFF$ обладает $\sS$-свойством (в терминах
потенциала~\eqref{2.6})
$$
\frac{\partial\bigl(3V^{\lambda}+G^{\lambda}_{\EE}+3g_S(\cdot,\infty)\bigr)}
{\partial n^{+}}(\zeta)
=
\frac{\partial\bigl(3V^{\lambda}+G^{\lambda}_{\EE}+3g_S(\cdot,\infty)\bigr)}
{\partial n^{-}}(\zeta),\quad\zeta\in \FF^0
$$
(для компакта $\EE$ соответствующее $\sS$-свойство вытекает из~\eqref{2.4.1}).
Соответствующим способом переформулируется
$\max$-$\min$-задача~\eqref{2.2}. При этом $\lambda=\lambda_F$.
\end{remark}

\begin{remark}\label{rem2}
Соотношение~\eqref{1.5} симметрично относительно $1$ и $f^2$ (т.е. обе его
части можно
поделить на $f^2$), поскольку функция $1/f$ принадлежит тому же классу что
и~$f$. Поэтому из~\eqref{1.5} сразу же вытекает, что слабая асимптотика для
$Q_{n,0}$ и $Q_{n,2}$ одна и та же. Для $Q_{n,1}$ ответ будет такой же, но
доказывать это придется отдельно.
\end{remark}

\begin{remark}\label{rem3}
Таким образом, в настоящей работе, в частности, утверждается, что точка
Чеботарёва {\it сохраняется неподвижной} при переходе от АП к Эрмита--Паде, а
компакт Шталя $S$ участвует в ответе (через функцию Грина $g_S(z,\infty)$).
Это означает, что и в общем случае произвольной алгебраической функции $f$
построить конденсатор Наттолла для набора из трех функций $1,f,f^2$
невозможно без знания структуры компакта Шталя для функции $f$. А
именно, конденсатор Наттолла строится с помощью
функции Грина $g_S(z,\infty)$ для области Шталя~$D$ и ``остова'' компакта
Шталя, который состоит из ``эффективных'' точек ветвления функции $f$
(см.~\cite{KovSu14}) и соответствующих им точек Чеботарёва.
Это утверждение вполне согласуется с результатами А.~И.~Аптекарева,
В.~А.~Калягина, В.~Г.~Лысова и В.~Н.~Тулякова
(см.~\cite{ApKa86}, \cite{ApKuVa08},
\cite{ApLy10},~\cite{ApLyTu11},~\cite{ApTu12b})
о {\it сильной асимптотике} полиномов Эрмита--Паде второго рода, традиционно
основанным на использовании кубической функции, введенной в~\cite{ApKa86}.
\end{remark}

\section{Вариационный метод}\label{s3}

\subsection{}\label{s3s1}
Основная цель настоящего параграфа -- исходя из п.1 утверждения~\ref{pro1}
(см.~\eqref{ext1}) вариационным методом (см.~\cite{PeRa94},
\cite{Apt08}, \cite{MaRa11a},~\cite{MaRa11b},
\cite{MaRaSu11a},~\cite{MaRaSu11b},~\cite{BuMaSu12},
\cite{Sue12a},~\cite{Bus13}) получить следующее
соотношение:
$$
\myFF=\biggl\{
w:\Re\int_{a_1}^w\sqrt{\frac{g(\zeta)}{A_4(\zeta)}}\,d\zeta=0
\biggr\},
$$
где $g$ -- некоторая функция, голоморфная в $\myo{\CC}\setminus{(S\cup{\EE})}$,
$A_4(\zeta)=(\zeta-v)\prod_{j=1}^3(\zeta-a_j)$.
Именно это соотношение приводит к п.2 утверждения~\ref{pro1}, а в дальнейшем
-- и к утверждениям~\ref{pro2} и~\ref{pro3}.

Рассмотрим следующую задачу. Пусть $\EE,\myFF$ -- допустимые компакты
($(\EE,\myFF)$ -- допустимый конденсатор), $\lambda=\lambda_{\EE,\myFF}\in M_1(\myFF)$ --
равновесная мера (с носителем на $\myFF$):
\begin{equation}
3V^\lambda(z)+G^\lambda_{\EE}(z)+3g_S(z,\infty)\equiv\gamma(\EE,\myFF),
\qquad z\in \myFF.
\label{3.1}
\end{equation}
Предположим, что допустимый конденсатор $(\EE,\myFF)$ является
{\it экстремальным}, т.е. что справедливо соотношение
\begin{equation}
J_\psi(\EE,\myFF,\lambda)
=\sup_{(K,K^{\perp}), }J_\psi(K,K^{\perp},\lambda_{K,K^{\perp}}),
\label{3.2}
\end{equation}
где $\lambda=\lambda_{E,F}$,
\begin{equation}
J_\psi(K,K^{\perp},\lambda)
:=\int_{K^{\perp}}\bigl(3V^\lambda(z)+G^\lambda_{\EE}(z)\bigr)\,d\lambda(z)
+6\int_{K^{\perp}}g_S(z,\infty)\,\lambda(z),
\label{3.3}
\end{equation}
здесь $\lambda=\lambda_{K,K^{\perp}}$, $\psi(z)=3g_S(z,\infty)$ -- внешнее поле.
Таким образом, (допустимый) конденсатор $(\EE,\myFF)$ --
{\it стационарная ``точка''} для функционала энергии~\eqref{3.3}, т.е. для
любой пары $(K,K^{\perp})$ выполняется неравенство:
$$
J_\psi(\EE,\myFF,\lambda_{\EE,\myFF})
\geq J_\psi(K,K^{\perp},\lambda_{K,K^{\perp}}).
$$

Пусть функция $h(z)$ -- произвольная функция, голоморфная на компакте
$\EE\cup \myFF$ и обращающаяся в ноль в точках $\{a_1,a_2,a_3,v\}$, $t\in\CC$ -- комплексный
параметр.
Введем следуя~\cite{PeRa94} (см. также~\cite{MaRa11a},~\cite{MaRaSu11b})
следующее преобразование (``вариацию'') $z\mapsto z_t=z+th(z)$. Зафиксируем
некоторую окрестность $U$ компакта $\EE\cup \myFF$ такую, что $h\in\HH(\myo{U})$.
Нетрудно увидеть, что при достаточно малых $t$, $0<|t|<\eps_0$, $\eps_0>0$,
преобразование $z\mapsto z_t$ однолистно в окрестности $U$. Преобразование
$z_t$ порождает естественное преобразование множеств, принадлежащих окрестности
$U$, $e\mapsto e_t$, и мер $\mu$, носители которых принадлежат $U$, по
правилу $\mu\mapsto\mu^t$, $\mu^t(e_t)=\mu(e)$.

\subsection{}\label{s3s2}
Выберем теперь
\begin{equation}
h(z)=\frac{A_4(z)}{(z-w)^6},
\label{var-h}
\end{equation}
где $A_4(z):=(z-v)\prod_{j=1}^3(z-a_j)$, а $w\in\CC\setminus{U}$ -- произвольная точка (параметр вариации).
Тогда для вариации логарифмического ядра и функции Грина имеем соответственно
(см.~\cite{MaRaSu11b},~\cite{Sue12a}):
\begin{gather}
\delta_t\log\frac1{|z-\zeta|}:=
\log\frac1{|z_t-\zeta_t|}-\log\frac1{|z-\zeta|}
=\Re\biggl\{\frac{t}{5!}\biggl(\frac{A_4(w)}{(z-w)(z-\zeta)}\biggr)^{(5)}_w
\biggr\}+O(t^2),
\label{3.4}\\
\delta_t g_{\EE}(z,\zeta):=g_{\EE_t}(z_t,\zeta_t)-g_{\EE}(z,\zeta)
=\Re\biggl\{\frac{t}{5!}\bigl(P'(w,z)P'(w,\zeta)A_4(w)\bigr)^{(5)}_w\biggr\}
+O(t^2),
\label{3.5}
\end{gather}
где $P(w,z)=P_E(w,z)=g_{\EE}(w,z)+i\myt g_{\EE}(w,z)$ -- (многозначная) комплексная
функция Грина для дополнения к компакту $\EE$, производная $P'(w,\cdot)$ берется
по переменному $w$, запись $(\cdot)^{(5)}_w$ здесь и всюду в дальнейшем означает
5-ую производную по $w$ выражения, стоящего в круглых скобках. Кроме того имеем
следующую формулу для разности
\begin{align}
\log\frac1{|z_t-\zeta|}-\log\frac1{|z-\zeta|}
&=-\Re\biggl\{\log\frac{z_t-\zeta}{z-\zeta}\biggr\}
=-\Re\biggl\{\log\biggl(1+\frac{th(z)}{z-\zeta}\biggr)\biggr\}\notag\\
&=-\Re\biggl\{\frac{th(z)}{z-\zeta}\biggr\}+O(t^2),
\label{3.6}
\end{align}
где как и выше функция $h$ задана формулой~\eqref{var-h}. Пусть теперь $\mu\in
M_1(\myFF)$ -- произвольная мера с носителем на~$\myFF$. Так как
$g_S(z,\infty)=\gamma_S-V^{\lambda_S}(z)$, то с помощью~\eqref{3.6}
получаем
\begin{align}
K(\mu):
&=\int_{\myFF_t}g_S(z,\infty)\,d\mu^t(z)-
\int_{\myFF}g_S(z,\infty)\,d\mu(z)\notag\\
&=\int_Sd\lambda_S(\zeta)\Re\biggl\{\int_{\myFF}\frac{th(z)}{z-\zeta}
\,d\mu(z)\biggr\}+O(t^2)\notag\\
&=\Re\biggl\{\int_{\myFF}\biggl\{\int_S\frac{d\lambda_S(\zeta)}{z-\zeta}
\biggr\}th(z)\,d\mu(z)\biggr\}+O(t^2)\notag\\
&=\Re\biggl\{t\int_{\myFF}\myh\lambda_S(z)h(z)\,d\mu(z)\biggr\}+O(t^2).
\label{3.6.1}
\end{align}
Так как
$$
d\lambda_S(\zeta)=\frac1{\pi i}\sqrt{\frac{\zeta-v}{A_3(\zeta)}}\,d\zeta,
$$
где $A_3(z)=\prod_{j=1}^3(z-a_j)$, то
\begin{equation}
\myh{\lambda}_S=\int_S\frac{d\lambda_S(\zeta)}{z-\zeta}
=\sqrt{\frac{z-v}{A_3(z)}}.
\label{3.6.2}
\end{equation}
Следовательно, для $h$, заданной формулой~\eqref{var-h}, в силу равенства
$$
A_4(z)\myh{\lambda}_S(z)=(z-v)^{3/2}\sqrt{A_3(z)},
$$
из~\eqref{3.6.1} и~\eqref{3.6.2} получаем
\begin{align}
K(\mu)
&=\Re\biggl\{t\int_{\myFF}\frac{A_4(z)}{(z-w)^6}
{\myh\lambda}_S(z)\,d\mu(z)\biggr\}+O(t^2)\notag\\
&=\Re\biggl\{\frac{t}{5!}\biggl(\int_{\myFF}\frac{A_4(z)}{z-w}
{\myh\lambda}_S(z)\,d\mu(z)\biggr)^{(5)}_w\biggr\}+O(t^2).
\label{3.6.3}
\end{align}
Теперь с учетом соотношений
\eqref{3.4}--\eqref{3.6} и~\eqref{3.6.3} для  произвольной меры~$\mu\in
M_1(\myFF)$ стандартным образом (ср.~\cite{PeRa94},
\cite[формула~(39)]{MaRaSu11b}) получаем
\begin{align}
\delta_t J_\psi(\EE,\myFF;\mu):
&=J_\psi(\EE_t,\myFF_t;\mu^t)-J_\psi(\EE,\myFF;\mu)\notag\\
&=\Re\biggl\{\frac{t}{5!}\biggl(\Bigl[3\myh{\mu}(w)^2+\sP'_{\EE,\mu}(w)^2\Bigr]
A_4(w)\biggr)^{(5)}_w\biggr\}+6K(\mu)+O(t^2),
\label{3.7}
\end{align}
где
$$
\myh\mu(z):=\int_{\myFF}\frac{d\mu(\zeta)}{z-\zeta},\qquad
\sP_{\EE,\mu}(z):=\int_{\myFF}P_{\EE}(z,\zeta)\,\mu(\zeta),
$$
$P_{\EE}(z,\zeta)=g_{\EE}(z,\zeta)+i\myt g_{\EE}(z,\zeta)$ -- многозначная функция Грина
для дополнения к компакту $\EE$, а для $K(\mu)$ имеет место
представление~\eqref{3.6.3}.

\subsection{}\label{s3s3}
Введем семейство мер $\sigma_t\in M_1(\myFF)$ таких, что
$(\sigma_t)^t=\lambda_{\myFF_t}$, и семейство мер
$\lambda^t:=(\lambda_{\myFF})^t$.
Тогда пользуясь произвольностью малого комплексного параметра~$t$ с помощью
этих 4-х мер из соотношений~\eqref{3.6.3} и~\eqref{3.7} для потенциала
равновесной меры~$\lambda=\lambda_{\EE,\myFF}$ стандартным образом
(см.~\cite{PeRa94}, \cite{MaRaSu11b}) получаем тождество:
\begin{equation}
\Bigl\{3\myh\lambda(w)^2+\sP'_{\EE,\lambda}(w)^2\Bigr\}A_4(w)
+6\int_{\myFF}\frac{A_4(z)\myh\lambda_S(z)}{z-w}\,d\lambda(z)
\equiv b_4(w),
\label{3.8}
\end{equation}
где $b_4(w)$ -- некоторый полином по $w$ 4-й степени. С помощью тождества
$A_4(z)=(A_4(z)-A_4(w))+A_4(w)$ и равенства
\begin{align}
\int_{\myFF}\frac{A_4(z)\myh\lambda_S(z)}{z-w}\,d\lambda(z)
&=\int_{\myFF}\frac{A_4(z)-A_4(w)}{z-w}\myh\lambda_S(z)\,d\lambda(z)
-\bigl(\myh\lambda_S(\cdot)\lambda\bigr)^{\myh{\ }}(w)A_4(w)\notag\\
&=-\bigl(\myh\lambda_S(\cdot)\lambda\bigr)^{\myh{\ }}(w)A_4(w)+c_4(w),
\end{align}
где $\bigl(\myh\lambda_S(\cdot)\lambda\bigr)^{\myh{\ }}$ -- преобразование
Коши соответствующей меры, а $c_4(w)$ -- полином 4-й степени по $w$,
соотношение~\eqref{3.8} приводится к виду
\begin{equation}
3\myh\lambda(w)^2+\sP_{\EE,\lambda}'(w)^2
-6\bigl(\myh\lambda_S(\cdot)\lambda\bigr)^{\myh{\ }}(w)
\equiv\frac{B_4(w)}{A_4(w)},
\label{3.9}
\end{equation}
где $B_4(w)$ -- некоторый полином 4-й степени,
$A_4(w)=(w-v)\prod_{j=1}^3(w-a_j)$.
Воспользуемся теперь тем, что для гринова потенциала меры $\lambda$ имеем
\begin{equation}
G^\lambda_{\EE}(z):=\int_{\myFF} g_{\EE}(z,\zeta)\,d\lambda(\zeta)
=V^\lambda(z)-V^{\beta(\lambda)}(z)+\const,
\label{3.10}
\end{equation}
где $\beta(\lambda)=\beta_{\EE}(\lambda)\in M_1(\EE)$ -- выметание меры $\lambda$
с компакта $\myFF$ на (вторую пластину конденсатора) $\EE$. Из~\eqref{3.10}
получаем
\begin{equation}
\sP_{\EE,\lambda}'(w)=\sV'_\lambda(w)-\sV'_{\beta(\lambda)}(w)
=\myh{\beta(\lambda)}(w)-\myh\lambda(w),
\label{3.11}
\end{equation}
где
$$
\sV_\mu(z):=-\int\log(z-\zeta)\,d\mu(\zeta).
$$
Кроме того имеем
\begin{equation}
\bigl(\myh\lambda_S(\cdot)\lambda\bigr)^{\myh{\ }}(w)
=\myh\lambda_S(w)\myh\lambda(w)
-\int_{\myFF}\frac{\myh\lambda_S(z)-\myh\lambda_S(w)}{z-w}\,d\lambda(z).
\label{3.12}
\end{equation}
Из~\eqref{3.9},~\eqref{3.11} и~\eqref{3.12} получаем следующее квадратное
уравнение на функцию $\myh\lambda(w)$:
\begin{equation}
4\myh\lambda(w)^2-2\myh{\beta(\lambda)}(w)\myh\lambda(w)
-6\myh\lambda_S(w)\myh\lambda(w)\equiv\frac{g_1(w)}{A_4(w)},
\label{3.13}
\end{equation}
где $g_1(w)=\myh{\beta(\lambda)}(w)^2$ -- функция, голоморфная вне $\EE$.
Разрешая~\eqref{3.13} относительно функции $\myh\lambda(w)$, получаем
соотношение
\begin{equation}
4\myh\lambda(w)-\myh{\beta(\lambda)}(w)-3\myh\lambda_S(w)
\equiv\sqrt{\frac{g(w)}{A_4(w)}},
\label{3.14}
\end{equation}
где $g(w)$ -- некоторая функция, голоморфная вне $\EE$ и $S$ и выбрана
надлежащая ветвь корня. Проинтегрируем теперь это соотношение~\eqref{3.14},
возьмем вещественную часть от обеих частей и воспользуется
определением (см.~\eqref{3.11}) выметания $\beta(\lambda)=\beta_{\EE}(\lambda)$
меры $\lambda$.
Тогда из~\eqref{3.14} получаем
\begin{equation}
3V^\lambda(w)+G^\lambda_{\EE}(w)-3V^{\lambda_S}(w)
=\Re\int_{a_1}^w\sqrt{\frac{g(\zeta)}{A_4(\zeta)}}\,d\zeta+\const.
\label{3.15}
\end{equation}
Наконец, используя тождеством $V^{\lambda_S}(z)=\gamma_{\Rob}-g_S(z,\infty)$
из~\eqref{3.15} получаем
\begin{equation}
3V^\lambda(w)+G^\lambda_{\EE}(w)+3g_S(w,\infty)
=\Re\int_{a_1}^w\sqrt{\frac{g(\zeta)}{A_4(\zeta)}}\,d\zeta+\const.
\label{3.16}
\end{equation}
В силу определения равновесной меры $\lambda=\lambda_{\myFF}$ имеем:
\begin{equation}
3V^\lambda(w)+G^\lambda_{\EE}(w)+3g_S(w,\infty)
\equiv\gamma(\EE,\myFF)\quad\text{на}\quad
\myFF.
\label{3.17}
\end{equation}
Из~\eqref{3.16} и~\eqref{3.17} вытекает, что $\const=\gamma(\EE,\myFF)$ и
\begin{equation}
\myFF=\biggl\{
w:\Re\int_{a_1}^w\sqrt{\frac{g(\zeta)}{A_4(\zeta)}}\,d\zeta=0\biggr\}.
\label{3.18}
\end{equation}
Непосредственно из~\eqref{3.18} вытекает $\sS$-свойство компакта $\myFF$.

\section{Обоснование утверждений~\ref{pro2} и~\ref{pro3}}\label{s4}

\subsection{}\label{s4s1}
Основная цель настоящего параграфа -- доказательство нижеследующего
Предложения~\ref{Prop1}, непосредственно из которого вытекает принципиальная
возможность использования общего метода Гончара--Рахманова~\cite{GoRa87} для
обоснования утверждений~\ref{pro2} и~\ref{pro3}.

Для дальнейшего договоримся о следующих соглашениях. Будем считать, что точки $a_1,a_2,a_3$
находятся в ``общем положении''. В частности, это означает, что для диагональных
АП $[n/n]_f$ функции~$f$ все достаточно далекие индексы {\it нормальны}.
Мы будем предполагать, что пластины конденсатора $(\EE,\myFF)$ не содержат
бесконечно удаленной точки (этого всегда можно добиться дробно-линейным
преобразованием комплексной плоскости).
Под
$$
\oint_S(\cdot),\qquad\oint_{\EE}(\cdot),\quad\oint_{\myFF}(\cdot),
$$
будем понимать интеграл от выражения в скобках по произвольному
контуру~$\gamma$, охватывающему соответственно компакт $S$, $\EE$ или~$\myFF$.
Для произвольной функции $f$ вида~\eqref{01.1} ее скачок на $S$ ($\EE$ или $\myFF$) может
иметь в концевых точках~$a_j$ неинтегрируемые особенности. После умножения
на некоторый полином фиксированной степени скачок будет уже интегрируемой
функцией. Поскольку мы изучаем слабую асимптотику, такое умножение не повлияет
на результаты работы, поэтому всюду в дальнейшем мы предполагаем интегрируемость
скачков.
Наконец, под $\mint(\cdot)$ будем понимать открытую составляющую соответствующей
дуги (без точек $\{a_1,a_2,a_3,v_1,v_2\}$; такое обозначение мы будем
использовать наряду с обозначением $(\cdot)^{0}$).

Непосредственно из~\eqref{1.5} вытекает, что выполняются следующие ``соотношения
ортогональности''
\begin{equation}
\int_\gamma(Q_{n,1}f+Q_{n,2}f^2)(\zeta)\cdot\zeta^k\,d\zeta=0,
\qquad k=0,1,\dots,2n,
\label{4.1}
\end{equation}
где $\gamma$ -- произвольный контур, отделяющий точки $a_j,v$ от бесконечно
удаленной точки. Отметим, что число свободных параметров в
соотношении~\eqref{4.1} равно $2n+1$ ($=$ суммарное число коэффициентов
полиномов $Q_{n,1}$ и $Q_{n,2}$ минус $1$). Тем самым с учетом определенной
нормировки соотношениями~\eqref{4.1} полиномы $Q_{n,1}$ и $Q_{n,2}$
определяются однозначно; полином $Q_{n,0}$ уже однозначно находится
из~\eqref{1.5}. Соотношения~\eqref{4.1} эквивалентны соотношениям
\begin{equation}
\int_\gamma(Q_{n,1}f+Q_{n,2}f^2)(\zeta)\cdot q(\zeta)\,d\zeta=0,
\label{4.2}
\end{equation}
где $q\in\PP_{2n}=\CC_{2n}[\zeta]$ -- произвольный полином степени $\leq{2n}$.
Избавимся теперь в~\eqref{4.2} от первого слагаемого,
содержащего функцию~$f$, под знаком интеграла в~\eqref{4.2}. Для этого в
качестве $q$ выберем полиномы Паде $P_{n+j,1}$ из~\eqref{01.7}, т.е. знаменатели
классических АП, соответствующих индексу $j=1,2,\dots,n$: $q=P_{n+j,1}$,
$j=1,2,\dots,n$. Поскольку $\mdeg{P_{n+j,1}}=n+j$ (напомним, что по
предположению все достаточно далекие индексы нормальны для диагонали таблицы
Паде), то в дальнейшем мы можем считать, что $P_{n+j,1}(z)=z^{n+j}+\dotsb$,
т.е. все $P_{n+j,1}$ -- монические полиномы соответствующей степени. В
дальнейшем для краткости будем вместо $P_{n+j,1}$ использовать обозначение
$P_{n+j}$.
Непосредственно из определения АП (см.~\eqref{01.7}) вытекает, что для любого
полинома $p\in\PP_{n}$ выполняются условия ортогональности
\begin{equation}
\int_\gamma P_{n+j,1}(\zeta)\cdot p(\zeta)f(\zeta)\,d\zeta=0,
\qquad j=1,2,\dots,n.
\label{4.3}
\end{equation}
Из~\eqref{4.2} и~\eqref{4.3} получаем следующие
{\it $n$ соотношений ортогональности} (ср.~\cite{Sue97}):
\begin{equation}
0=\int_\gamma(Q_{n,2}f^2)(\zeta)\cdot P_{n+j,1}(\zeta)\,d\zeta
=\oint_{\EE}(Q_{n,2}f^2)(\zeta)\cdot P_{n+j,1}(\zeta)\,d\zeta,
\qquad j=1,2,\dots,n.
\label{4.4}
\end{equation}

Справедлива следующее

\begin{Proposi}\label{Prop1}
Соотношения ортогональности~\eqref{4.4},
$$
\oint_{\EE} Q_{n,2}(\zeta)\cdot P_{n+j,1}(\zeta)\cdot f^2(\zeta)\,d\zeta=0,
\qquad j=1,2,\dots,n,
$$
эквивалентны следующему соотношению
\begin{equation}
0=\int_{\myFF}Q_{n,2}(t)\cdot\frac{\omega_n(t)}{q_m(t)}
\biggl\{\int_{\EE}\frac{q^2_m(u)\,d\rho_n(u)}{\omega_n(u)(t-u)}
\biggr\}
R_{n+m}(t)\cdot f^{*}(t)\,dt,
\label{GoRa1}
\end{equation}
где $\omega_n$ -- произвольный полином степени $\leq{n-2m_0}$ ($m_0$ --
фиксированное число) с нулями на произвольном зафиксированном компакте~$K$,
$K\subset\myo{\CC}\setminus({S}\cup{E})$, комплексная мера $d\rho_n(\zeta)$
имеет вид
$$
d\rho_n(\zeta)=\frac{f_j(\zeta)\,d\zeta}
{r^{+}_n(\zeta)r^{-}_n(\zeta)},\qquad \zeta\in\EE^0_j,
\quad j=1,2,3,
$$
где $f_j(\zeta)=(f^{+}+f^{-})(\zeta)$, $\zeta\in \EE_j$, $j=1,2,3$,
$r_n$ -- нормированная функция остатка, функция $f^{*}(t)$ определена
на~$\myFF$ (см.~ниже~\eqref{f-star}), а $q_m$ -- полином степени $\leq n/2-m_0$,
ортогональный на пластине $\EE$ относительно меры
$d\rho_n(\zeta)/\omega_n(\zeta)$. При этом
$$
d\rho_n(\zeta)\to f_j(\zeta)\,d\zeta,\quad n\to\infty,\quad \zeta\in\EE^0_j,
\quad j=1,2,3.
$$
\end{Proposi}

\begin{proof}
Так как в~\eqref{4.4} $\gamma$ -- произвольный контур, охватывающий точки
$a_1,a_2,a_3,v$, то эти соотношения~\eqref{4.4} при тех же $j=1,2,\dots,n$
эквивалентным образом переписываются заменой интегрирования по $\gamma$ на
интегрирование по компакту Шталя~$S$ следующим образом:
\begin{equation}
\int_{S^{+}}Q_{n,2}(\zeta)P_{n+j,1}(\zeta)(f^{+}+f^{-})(\zeta)\,d\sigma(\zeta)=0,
\qquad j=1,2,\dots,n,
\label{4.5}
\end{equation}
где $d\sigma(\zeta):=(f^{+}-f^{-})(\zeta)\,d\zeta=\Delta{f}(\zeta)\,d\zeta$,
$\zeta\in{S}$, на дугах $S_j=(a_j,v)$ компакта $S$ выбрана положительная
ориентация ``от точки ветвления $a_j$ к точке~$v$''.\footnote{Так же как
выбрано в работе~\cite{MaRaSu12}.} Тем самым задана
и ориентация $S^{+}$ всего компакта Шталя. Под ``положительной'' частью
прилегающей к $S_j$ части области $D$ в окрестности $S_j$ понимается та
``половина'', которая остается ``слева'' при движении по $S_j$ в положительном
направлении ``от точки $a_j$ к точке~$v$'', вторая ``половина'' (которая
остается справа) считается отрицательной частью. Соответственно этому для точек
$\zeta\in{S^0_j}$ под $f^{+}(\zeta)$ понимаются предельные значения функции
$f(z)$ при $z\to\zeta$ из ``положительной'' части, а под $f^{-}(\zeta)$ --
предельные значения функции $f(z)$ при $z\to\zeta$ из ``отрицательной'' части
прилегающей части области~$D$. Отметим, что скачок
$\Delta{f}(\zeta)=(f^{+}-f^{-})(\zeta)$, $\zeta\in S^0$, -- (кусочно)
голоморфная функция на каждой открытой дуге $S^0_j$, $j=1,2,3$; то же самое
справедливо и для функции $f_j(\zeta)$, $\zeta\in{S^0_j}$, $j=1,2,3$, где
\begin{equation}
f_j(\zeta):=(f^{+}+f^{-})(\zeta)\bigr|_{S^{0}_j},\qquad
j=1,2,3.
\label{4.6}
\end{equation}

Пусть $D=D_1\sqcup D_2\sqcup D_3$ -- произвольное открытое множество, состоящее
из трех неналегающих областей, $D_j\cap D_k=\varnothing$ при $k\neq{j}$, каждая
из которых содержит (``нанизана'' на) соответствующую (``свою'') дугу Шталя
$S^{0}_j$, $D_j\supset
S^{0}_j$, тем самым, $\partial D_j\ni\{a_j,v\}$, $j=1,2,3$. Очевидно, что
функция $f_j$ голоморфно продолжается с дуги $S^{0}_j$ в соответствующую ей
область $D_j$. Предположим, что области $D_j$ имеют аналитические границы
$\gamma_j=\partial{D_j}$, ориентированные стандартным образом: при их обходе в
положительном направлении $\gamma^{+}_j$ соответствующая область $D_j$ остается
слева. Пусть
\begin{equation}
R_n(z)=R_n(z;\sigma):=\int_{S^{+}}\frac{P_n(\zeta;\sigma)\,d\sigma(\zeta)}{z-\zeta},
\qquad z\notin{S},
\label{4.7}
\end{equation}
-- функция второго рода, соответствующая ортогональному относительно веса
$d\sigma(\zeta)=\Delta f_j(\zeta)\,d\zeta$, $\zeta\in S^0_j$, полиному
$P_{n,1}(\zeta)=P_n(\zeta;\sigma)$ (см.~\eqref{4.5}, иначе говоря,
$R_n(z)$ -- функция остатка для диагональной АП $[n/n]_f(z)$).
Отметим, что функция $R_n(z)$ имеет в бесконечно удаленной точке нуль
порядка\footnote{Напомним, что мы предполагаем все достаточно далекие индексы
нормальными.}~$n+1$:
\begin{equation}
R_n(z)=O\(\frac1{z^{n+1}}\),\qquad z\to\infty.
\label{4.7.1}
\end{equation}

Нетрудно увидеть, что в терминах функций $R_n$ соотношения~\eqref{4.5}
переписываются в следующем виде:
\begin{equation}
\sum_{j=1}^3\int_{\gamma_j}Q_{n,2}(t)R_{n+k}(t)f_j(t)\,dt=0,
\qquad k=1,\dots,n,
\label{4.8}
\end{equation}
где $\gamma_j=\partial D_j$.
Для ортогональных полиномов $P_{n}(z)=z^n+\dotsb$ справедливы следующие трехчленные
рекуррентные соотношения:
\begin{equation}
P_n(z)=(z-b_n)P_{n-1}(z)-a^2_{n-1}P_{n-2}(z),
\qquad n\geq n_0,
\label{4.9}
\end{equation}
где все $a_n\neq0$.
Из~\eqref{4.7} вытекает, что такие же соотношение, но с другими начальными
данными справедливы и для функций второго рода $R_n(z)$:
\begin{equation}
R_n(z)=(z-b_n)R_{n-1}(z)-a^2_{n-1}R_{n-2}(z).
\label{4.10}
\end{equation}
Перепишем теперь соотношения~\eqref{4.8} в эквивалентном виде
\begin{equation}
\sum_{j=1}^3\int_{\gamma_j}Q_{n,2}(t)
\biggl(\sum_{k=1}^n c_kR_{n+k}(t)\biggr)f_j(t)\,dt=0,
\label{4.11}
\end{equation}
где $c_k\in\CC$, $k=1,\dots,n$, -- произвольные числа. Теперь с помощью
рекуррентных соотношений~\eqref{4.10} выполним следующее преобразование:
\begin{equation}
\sum_{k=1}^n c_kR_{n+k}(t)=p_m(t)R_{n+m}(t)+q_{m}(t)R_{n+m+1}(t),
\label{4.12}
\end{equation}
где, например, для четного $n$ индекс $m=n/2$, а полиномы
$p_m,q_m\in\PP_{m-1}$. Отметим, что поскольку в~\eqref{4.11} постоянные $c_k$
-- произвольные, то и полиномы $p_m,q_m$ --
{\it произвольные} полиномы с комплексными коэффициентами, $p_m,q_m\in\PP_{m-1}$.
С учетом равенства~\eqref{4.12} соотношение~\eqref{4.11} принимает вид:
\begin{equation}
0=
\sum_{j=1}^3\int_{\gamma^{+}_j}Q_{n,2}(t)
\biggl\{p_m(t)R_{n+m}(t)+q_{m}(t)R_{n+m+1}(t)\biggr\}f_j(t)\,dt.
\label{4.13}
\end{equation}
Заметим теперь, что в силу~\eqref{4.7.1} при любых $p_m,q_m\in\PP_{m-1}$
выражение в фигурных скобках в~\eqref{4.13} имеет в бесконечно удаленной
точке нуль порядка $n+2$, а оставшаяся часть подынтегрального выражения --
полюс порядка $\leq{n}$. Следовательно, интегралы в~\eqref{4.13} не
изменятся, если мы продеформируем ориентированные контуры $\gamma^{+}_j=\partial D_j$ в
соответствующим образом ориентированные контуры $\partial D^{*}_j$, где области
$D^{*}_1\sqcup D^{*}_2\sqcup D^{*}_3$  -- дают решение экстремальной
задачи~\eqref{2.2}. Теперь уже (с учетом ориентации) интегрирование по
$\gamma_j$ в~\eqref{4.13} можно заменить на интегрирование по состоящим из аналитических
дуг экстремальным компактам $\myFF_k$, $k=1,2,3$. Точнее, определим функцию
$f^{*}_k(t)$, $t\in \myFF_k$, $k=1,2,3$, как разность значений (``скачок'') функций
$f_{j_1}(\zeta_1)$ и $f_{j_2}(\zeta_2)$, заданных в двух разных областях
$D_{j_1}\ni\zeta_1$ и $D_{j_2}\ni\zeta_2$, $j_1\neq j_2$, при
$\zeta_1,\zeta_2\to t\in \myFF_k$. Нетрудно видеть, что при таком определении для
произвольной функции $f$ вида~\eqref{01.1} имеем: $f^{*}_k$ -- аналитична на
$\FF^0_k$
и $f^{*}_k(t)\not\equiv0$ при $t\in\FF^0_k$, $k=1,2,3$. Таким образом,
соотношение~\eqref{4.13} преобразуется следующим образом:
\begin{equation}
0=
\sum_{j=1}^3\int_{\myFF_k}Q_{n,2}(t)
\biggl\{p_m(t)R_{n+m}(t)+q_{m}(t)R_{n+m+1}(t)\biggr\}f^{*}_k(t)\,dt.
\label{4.14}
\end{equation}
Отметим, что ориентация в~\eqref{4.14} согласована с ориентацией границ
$\partial D^{*}_j$ и определением функций $f^{*}_k$, $k=1,2,3$, но для
дальнейших рассуждений это уже не существенно. Нам важно лишь то, что для
всех $k=1,2,3$ имеем: функция $f^{*}_k(t)\not\equiv0$ и аналитична на
$\FF^0_k$.

Введем функцию $f^{*}(t)$:
\begin{equation}
f^{*}(t)=f^{*}_k(t), \quad t\in{\myFF^0_k},\quad k=1,2,3,
\label{f-star}
\end{equation}
и перепишем соотношение~\eqref{4.14} в следующем эквивалентном виде:
\begin{equation}
0=
\int_{\myFF}Q_{n,2}(t)
\biggl\{p_m(t)R_{n+m}(t)+q_{m}(t)R_{n+m+1}(t)\biggr\}f^{*}(t)\,dt,
\label{4.15}
\end{equation}
где $\myFF$ -- (второй) экстремальный компакт в смысле задачи~\eqref{2.2}.
Подчеркнем, что в~\eqref{4.15} $p_m,q_m$ -- {\it произвольные}
полиномы с комплексными коэффициентами, $p_m,q_m\in\PP_{m-1}$.
Перепишем теперь соотношение~\eqref{4.15} еще раз в следующем эквивалентном
виде:
\begin{equation}
0=
\int_{\myFF}Q_{n,2}(t)
\biggl\{q_{m}(t)\frac{R_{n+m+1}}{R_{n+m}}(t)-p_m(t)\biggr\}R_{n+m}(t)f^{*}(t)
\,dt
\label{4.15.2}
\end{equation}
(при переходе от~\eqref{4.15} к~\eqref{4.15.2} поменяли знак у {\it произвольного} полинома
$p_m\in\PP_{m-1}$).


Непосредственно из определения~\eqref{01.7} функции остатка $R_n$ вытекает, что
это --
многозначная аналитическая функция на римановой сфере $\myo{\CC}$ с проколами
в точках $a_,a_2,a_3$; точнее, $R_n$ допускает многозначное аналитическое
продолжение по любому пути в $\myo{\CC}$, не пересекающему точки
$\{a_1,a_2,a_3\}$. При этом в области Шталя $D_S=\myo\CC\setminus{S}$ условием
$R_n(z)=O(1/z^{n+1})$,
$z\to\infty$, однозначно выделяется голоморфная (однозначная аналитическая)
ветвь функции $R_n$,  имеющая в бесконечно удаленной точке нуль кратности
в точности $n+1$ и возможно еще один ``свободный'' нуль в точке $\myt
z_n$, динамика которого при $n\to\infty$ подчиняется проблеме обращения Якоби
(см.~\cite[формула~(4.56)]{Nut86}). Тем самым функция $\pfi_n(z):=R_{n+1}(z)/R_n(z)$ --
мероморфная функция в области Шталя, с единственным возможным полюсом в точке
$z=\myt z_n$; обозначим через $c_n$ вычет этой функции в точке $z=\myt z_n$, допуская
при этом возможность $c_n=0$.
Пусть $\gamma^{+}=\gamma^{+}_1\cup\gamma^{+}_2\cup\gamma^{+}_3$ --
произвольная допустимая (ориентированная) кривая из класса $\GG$. Тогда функция
$\pfi_n(z)$ продолжается как (однозначная) мероморфная функция из окрестности
точки $z=\infty$ в область $\myo{\CC}\setminus\gamma$ и на двух разных берегах
$\gamma^{+}$ и $\gamma^{-}$ определены предельные значения $\pfi^{+}_n(z)$ и
$\pfi^{-}_n(z)$,
$z\in\gamma$, функции $\pfi_n(z)$ (отметим, что эти предельные функции --
аналитические на $\gamma^0$).

Пусть $r_n(z)=\alpha_n R_n(z)$ -- нормированные функции остатка, соответствующие
ортонормированным многочленам $p_{n,1}(z)=\alpha_n P_{n,1}(z)$, $r^{+}_n(z)$
и $r^{-}_n(z)$ -- соответствующие предельные значения функции $r_n(z)$ при
$z\to\zeta^{+}\in\gamma^{+}$ и
$z\to\zeta^{-}\in\gamma^{-}$.

Покажем теперь, что
отношение $R_{n+1}(z)/R_n(z)$ функций второго рода является ``почти''
марковской функцией. А именно, докажем справедливость следующего утверждения.

\begin{lemma}\label{le1}
Пусть функция $f$ вида~\eqref{01.1},
$\gamma=\gamma_1\cup\gamma_2\cup\gamma_3\in\GG$ -- произвольная кривая из
класса $\GG$. Тогда для $z\in\myo{\CC}\setminus\gamma$ при любом $n\in\NN$ для
отношения $R_{n+1}/R_n$ функций второго рода справедливо следующее
представление
\begin{equation}
\frac{R_{n+1}}{R_n}(z)
=\int_{\gamma}\frac{d\rho_n(\zeta)}{z-\zeta}+\frac{c_n}{z-\myt z_n}
=\mhat{\rho}_n(z)+\frac{c_n}{z-\myt z_n},
\label{5.2}
\end{equation}
где комплексная мера $d\rho_n(\zeta)$ имеет вид
\begin{equation}
d\rho_n(\zeta)=\frac{f_j(\zeta)\,d\zeta}
{r^{+}_n(\zeta)r^{-}_n(\zeta)},\qquad \zeta\in\gamma^0_j,
\quad j=1,2,3,
\label{5.3}
\end{equation}
при этом
$$
d\rho_n(\zeta)\to f_j(\zeta)\,d\zeta,
\quad n\to\infty,\quad \zeta\in\EE^0_j,
\quad j=1,2,3.
$$
\end{lemma}

\begin{proof}
Для $z\in\gamma$, $\gamma\in\GG$, имеем
\begin{equation}
\Delta_n(z):=\frac{R_{n+1}^{+}(z)}{R_{n}^{+}(z)}
-\frac{R_{n+1}^{-}(z)}{R_{n}^{-}(z)}
=\frac{R_{n+1}^{+}(z)R_{n}^{-}(z)-R_{n+1}^{-}(z)R_{n}^{+}(z)}
{R_n^{+}(z)R_{n}^{-}(z)}.
\label{5.4}
\end{equation}
Воспользуемся теперь рекуррентными соотношениями~\eqref{4.10}, которые запишем
для разных переменных $z$ и $\zeta$:
\begin{equation}
\begin{aligned}
R_{n+1}(z)&=(z-b_{n+1})R_n(z)-a^2_n R_{n-1}(z),\\
R_{n+1}(\zeta)&=(\zeta-b_{n+1})R_n(\zeta)-a^2_n R_{n-1}(\zeta).
\end{aligned}
\label{5.5}
\end{equation}
Из пары тождеств~\eqref{5.5} получаем следующее равенство
$$
R_{n+1}(z)R_n(\zeta)-R_{n+1}(\zeta)R_{n}(z)
=(z-\zeta)R_n(z)R_n(\zeta)-a^2_n[R_n(z)R_{n-1}(\zeta)-R_n(\zeta)R_{n-1}(z)],
$$
из которого вытекает, что для граничных значений функций остатка имеет место
следующее рекуррентное соотношение:
\begin{equation}
R_{n+1}^{+}(z)R_{n}^{-}(z)-R_{n+1}^{-}(z)R_{n}^{+}(z)
=a^2_n[R_{n}^{+}(z)R_{n-1}^{-}(z)-R_{n}^{-}(z)R_{n-1}^{+}(z)].
\label{5.6}
\end{equation}
С учетом начальных условий $R_{-1}(z)\equiv1$, $R_0(z)=\myh\sigma(z)$ и равенств
$R^{+}_0(z)-R^{-}(z)=-2\pi i\sigma'(z)=-2\pi if_j(z)$, $z\in\gamma_j$,
$j=1,2,3$, и
$a_k=\alpha_{k-1}/\alpha_k$ из~\eqref{5.6} получаем, что
\begin{equation}
R_{n+1}^{+}(z)R_{n}^{-}(z)-R_{n+1}^{-}(z)R_{n}^{+}(z)
=\prod_{k=1}^n a^2_k\cdot(-2\pi i\sigma'(z))
=\frac{-2\pi i f_j(z)}{\alpha^2_n},\qquad z\in\gamma_j.
\label{5.7}
\end{equation}
Переходя к нормированным функциям остатка $r_n(z)$ из~\eqref{5.4}
и~\eqref{5.7} получаем тождество
\begin{equation}
\Delta_n(z):=\frac{R_{n+1}^{+}(z)}{R_{n}^{+}(z)}
-\frac{R_{n+1}^{-}(z)}{R_{n}^{-}(z)}
=\frac{-2\pi i f_j(z)}{\alpha^2_n R_n^{+}(z)R_{n}^{-}(z)}
=\frac{-2\pi i f_j(z)}{r_n^{+}(z)r_{n}^{-}(z)},
\qquad z\in\gamma_j,\quad j=1,2,4.
\label{5.8}
\end{equation}
Следовательно, по теореме Коши имеем
\begin{align}
\frac{R_{n+1}}{R_n}(z)-\frac{c_n}{z-\myt z_n}
&=\frac1{2\pi i}\oint_\gamma\frac{R_{n+1}}{R_n}(\zeta)\frac{d\zeta}{\zeta-z}
=\frac1{2\pi i}\int_\gamma\Delta_n(\zeta)\frac{d\zeta}{\zeta-z}\notag\\
&=\int_\gamma\frac{d\rho_n(\zeta)}{z-\zeta}=\myh\rho_n(z),\quad z\notin\gamma,
\label{5.9}
\end{align}
где
\begin{equation}
d\rho_n(\zeta)
=\frac{f_j(\zeta)\,d\zeta}{r_n^{+}(\zeta)r_{n}^{-}(\zeta)},\qquad\zeta\in\gamma_j,
\quad j=1,2,4.
\label{5.10}
\end{equation}
Лемма~\ref{le1} доказана.
\end{proof}

\begin{remark}\label{rr1}
Напомним, что для классического случая, когда $S=[-1,1]$ -- единичный отрезок
и рассматриваются полиномы Чебышёва (первого рода) $T_n(z)
=\Phi(z)^n+\Phi(z)^{-n}$, ортогональные на $[-1,1]$ ($\Phi(z)=z+\sqrt{z^2-1}$
-- функция, обратная функции Жуковского), соответствующая функция остатка
имеет следующий явный вид:
$$
R_n(z)=\frac1{\sqrt{z^2-1}}\Phi^{-n}(z).
$$
Следовательно, для отношения $R_{n+1}(z)/R_n(z)$ имеем:
$$
\frac{R_{n+1}}{R_n}(z)=\frac1{\Phi(z)}=z-\sqrt{z^2-1}
=\int_{-1}^1\frac{d\sigma(x)}{z-x},\quad z\notin S,
$$
где $d\sigma(x)=2\sqrt{1-x^2}\,dx$ -- не зависит от~$n$.
\end{remark}


Соотношение~\eqref{GoRa1} уже вполне пригодно для применения
{\it общего метода Гончара--Рахманова}, разработанного в~\cite{GoRa87} и
там же примененного для полного решения широко известной задачи Варги~``об~$1/9$''
(отметим, что на {\it неодносвязный} случай этот метод обобщен в работе
В.~И.~Буслаева~\cite{Bus13}). Соответствующим рассуждениям, основанным
на~\eqref{4.15.2}, посвящен пункт~\ref{s4s3} нашей работы.

\subsection{}\label{s4s2}
В этом пункте для полноты изложения мы
опишем кратко узловые моменты общего метода работы~\cite{GoRa87}, где при
доказательстве общей теоремы~2 из~\cite{GoRa87}, непосредственно приводящей к
решению известной задачи ``об~$1/9$'', впервые было введено и использовано понятие компакта, {\it
симметричного во внешнем поле} (отметим еще раз, что для случая, когда
дополнение к симметричному компакту не является областью, этот метод был
обобщен в работе~\cite{Bus13}).

Свойство симметрии стационарного компакта (т.е. стационарной ``точки'' для
функционала энергии типа~\eqref{2.1})
является в~\cite{GoRa87} ключевым при
доказательстве п.~(i) теоремы~3 (см.~\cite[\S~2, формулы~(28)--(31)]{GoRa87} и
при доказательстве сходимости многоточечных АП по емкости
(см.~\cite[\S~3, п.~8, формулы~(37)--(38)]{GoRa87}). Точнее, речь идет о
{\it естественных
дополнительных условиях} на последовательности полиномов $P_n$, весовых
функций $\Psi_n$ и кривые $L$, нужные для того, чтобы утверждение нижеследующей
леммы~7 из~\cite{GoRa87} оставалось в силе при замене интеграла
\begin{equation}
\int_L|P_n(t)\Psi_n(t)\chi(t)\,dt|\quad\text{на}\quad
\biggl|\int_L P_n(t)\Psi_n(t)\chi(t)\,dt\biggr|
\label{s1.01}
\end{equation}
($\chi$ -- суммируемая функция на кривой~$L$).

\begin{LemmaS}[из работы~\cite{GoRa87}]
Пусть $L$ есть объединение конечного числа спрямляемых кривых в $\CC$, функции
$\Psi_n$ непрерывны на $L$ и удовлетворяют условию
\begin{equation}
\psi_n:=\frac1n\log\frac1{|\Psi_n|}\rightrightarrows\psi,\qquad z\in L
\label{le7.1}
\end{equation}
(предполагается, что $\Psi_n\neq0$ на $L$), функция $\chi$ такова, что $|\chi|$
суммируем и положителен почти всюду на~$L$. Если последовательность многочленов
$P_n$ удовлетворяет условию $\nu_n=\dfrac1n\mu(P_n)\to\nu$, то справедливо
соотношение
\begin{equation}
\lim_{n\to\infty}\biggl(\int_L|(P_n\Psi_n\chi)(t)\,dt|\biggr)^{1/n}=e^{-m},
\qquad m=\min_{L}(V^\nu+\psi).
\label{le7.2}
\end{equation}
\end{LemmaS}

Для замены~\eqref{s1.01} необходимо перейти от произвольных спрямляемых кривых
к аналитическим дугам $L$ и наложить определенное условие ``симметрии''
на $P_n\Psi_n$ (относительно~$L$) в окрестности точки минимума предельного
потенциала $V^\nu+\psi$, где $\frac1n\mu(P_n)\to\nu$,
$-\frac1n\log|\Psi_n(z)|\to\psi(z)$, функция $\psi$ -- гармоническая в
некоторой окрестности кривой~$L$. Следующие определения позволяют
сформулировать соответствующее требование симметрии.

Пусть $L$ -- простая аналитическая дуга и $z_0$ -- внутренняя точка~$L$.
Существуют окрестность $\sU$ точки $z_0$ и конформное отображение $\pfi\colon \sU\to
\mDD$ ($\mDD$ -- открытый единичный круг), такое, что
$$
\pfi(\sU\cap L)=(-1,1),\quad \pfi(z_0)=0.
$$
Это позволяет определить антиконформное отображение $\sU$ на себя
(локальный аналог комплексного сопряжения): $z\mapsto z^{*}=\pfi^{-1}(\myo{\pfi(z)})$
(черта, как обычно, означает комплексное сопряжение). Окрестности $\sU$ такого
типа будем называть $*$-{\it симметричными}. Отображение $z\mapsto z^{*}$
позволяет стандартным образом ввести отображение мер $\mu\mapsto\mu^{*}$,
заданных в~$\sU$: $d\mu^{*}(z)=d\mu(z^{*})$.

Вопрос, поставленный в~\eqref{s1.01},
решается в работе~\cite{GoRa87} на основе следующей леммы~9; $\sS$-свойство компакта~$\EE$ в
поле~$\psi$ используется только в этой лемме. А именно, используется вытекающее
из $\sS$-свойства свойство симметрии потенциала $V^\lambda+\psi$ в окрестностях
точек из $\Gamma^0$ относительно соответствующих $*$-отображений.
Точнее, {\it пусть $(\EE,\psi)\in S$ и $\Gamma=\Supp{\lambda}$; тогда для каждой
точки $z_0\in\Gamma^{0}$ существует $*$-симметричная окрестность~$\sU$ такая,
что}
\begin{equation}
(V^\lambda+\psi)(z)=
(V^\lambda+\psi)(z^{*}),\qquad z\in\sU.
\label{6.0}
\end{equation}
Свойство~\eqref{6.0} позволяет получить слабую асимптотику
\begin{equation}
\lim_{n\to\infty}|I_n(z)|^{1/n}
\label{6.1}
\end{equation}
для интегралов вида
\begin{equation}
I_n(z):=\frac1{Q_n(z)p_n(z)}
\oint_{\EE}(Q_np_n)(t)\frac{\Psi_n(t)\chi(t)\,dt}{t-z},
\label{6.2}
\end{equation}
где $Q_n$ -- многочлен степени~$n$, ортогональный с переменным
(зависящим от номера $n$) весом, в качестве $p_n$ допускается
брать {\it произвольный} многочлен степени~$\leq{n}$
(формулу~\eqref{6.1} мы воспроизвели здесь из~\cite[п.~8]{GoRa87} вместе
с обозначениями, они никак не связаны с обозначениями настоящей работы).
С помощью $\sS$-свойства во внешнем поле в~\cite{GoRa87} показано, что для
интеграла~\eqref{6.2} предел~\eqref{6.1} существует. Это утверждение --
один из ключевых моментов работы~\cite{GoRa87}.

Следующая лемма из работы~\cite{GoRa87} является ключевой для результатов
работы~\cite{GoRa87}: только в ней используется $\sS$-свойство стационарного
компакта $\EE$ для логарифмического функционала энергии с гармоническим внешним
полем $\psi(z)$ и вытекающее из него соотношение~\eqref{6.0}.

\begin{LemmaN}[из работы~\cite{GoRa87}]
Пусть $(\EE,\psi)\in{S}$, мера $\mu$, $|\mu|\leq1$, отлична от равновесной
меры $\lambda=\lambda(\EE,\psi)$ и $e$ -- фиксированный компакт нулевой емкости.
Тогда существуют точка $z_0\in\Gamma_0\setminus{e}$, ее $*$-симметричная
окрестность $\sU$ и мера $\nu$ такие, что выполняются следующие свойства:

(i) $|\nu|<1$, $S_\nu\subset(\EE\cup\sU)$;

(ii) $\nu|_{\sU}=(\mu|_{\sU})^{*}$;

(iii) $V^\nu(z)\equiv-\infty$, $z\in B=e\cup A\cup A_1$;

(iv) $(V^{\mu+\nu}+2\psi)(z^{*})=(V^{\mu+\nu}+2\psi)(z)$, $z\in\sU$;

(v) $(V^{\mu+\nu}+2\psi)(z_0)<(V^{\mu+\nu}+2\psi)(z)$,
$z\in \EE\setminus\{z_0\}$.
\end{LemmaN}


\subsection{}\label{s4s3}
Вернемся к нашей задаче.

Будем говорить, что компакт $K$, $K\cap{S}=\varnothing$, является {\it
допустимым} для (кусочно) голоморфной функции
$f_j(\zeta)=(f^{+}+f^{-})(\zeta)$,
$\zeta\in S^{0}_j$, $j=1,2,3$, если $\myo{\CC}\setminus{K}$ состоит из трех
компонент $D_1,D_2,D_3$ таких, что $D_j\supset{S^{0}_j}$ и $f_j\in\sH(D_j)$,
$j=1,2,3$. Семейство всех допустимых компактов будем обозначать
через~$\kk^{*}_f$.

\begin{definition}\label{def2}
Будем говорить, что допустимый компакт $\Gamma\in\kk^{*}_f$ обладает
$\sS$-{\it свойством} (или является $\sS$-{\it компактом}) относительно
потенциала $3V^\mu+G_E^\mu+\psi$, $\psi(z)=3g_S(z,\infty)$, если выполняются
следующие условия: $\Gamma$ состоит из конечного числа аналитических дуг и для
равновесной меры $\lambda_\Gamma\in M_1(\Gamma)$ компакта $\Gamma$,
\begin{equation}
3V^{\lambda_\Gamma}(z)+G^{\lambda_\Gamma}_E(z)+\psi(z)\equiv\gamma_\Gamma,\qquad z\in\Gamma,
\label{a3}
\end{equation}
выполняется равенство
\begin{equation}
\frac{\partial (3V^{\lambda_\Gamma}+G^{\lambda_\Gamma}_E+\psi)}{\partial n^+}(z)=
\frac{\partial (3V^{\lambda_\Gamma}+G^{\lambda_\Gamma}_E+\psi)}{\partial n^-}(z),
\quad z\in\Gamma^{0},
\label{a4}
\end{equation}
где $\Gamma^{0}$ -- совокупность открытых дуг, составляющих $\Gamma$,
$\partial/\partial n^{\pm}$ -- нормальные производные, взятые с противоположных
сторон $\Gamma$; для такого компакта $\Gamma=\Gamma(f)$ будем писать
$\Gamma\in\sS=\sS(f)$.
\end{definition}

Соотношения~\eqref{a3} (для $K=\Gamma$) и~\eqref{a4}, выполненные на компакте
$\Gamma$,
позволяют стандартным образом (см.~\cite{GoRa87})
для произвольной точки $z_0\in\Gamma^0$ определить такое антиконформное
отображение $z\mapsto z^{*}$ некоторой окрестности $\sU$ этой точки на себя,
при котором
\begin{equation}
(3V^{\lambda_\Gamma}+G^{\lambda_\Gamma}_E+\psi)(z)=
(3V^{\lambda_\Gamma}+G^{\lambda_\Gamma}_E+\psi)(z^{*}),\qquad z\in\sU
\label{a1.3}
\end{equation}
(операция ``$*$'' -- локальный аналог комплексного сопряжения; см.
лемму~\ref{lemma2}).
Соотношение~\eqref{a1.3} -- непосредственное следствие
$\sS$-свойства компакта $\Gamma\in\kk^{*}_f$, также как и
в~\cite{GoRa87} именно оно лежит в основе рассуждений, приводящих к
утверждению~\ref{pro2}.

Вопрос, поставленный в~\eqref{s1.01}, решается на основе следующих
леммы~\ref{lemma2} и утверждения~\ref{pro4}; $\sS$-свойство компакта~$\Gamma$ в поле~$\psi$
используется только в этих леммах. А именно, используется вытекающее из
$\sS$-свойства~\eqref{a4} свойство симметрии потенциала
$3V^\lambda+G^\lambda_E+\psi$ в окрестностях
точек из $\mint F$ относительно соответствующих $*$-отображений\footnote{В
отличие от работы~\cite{GoRa87}, здесь мы имеем дело только с такой ситуацией,
когда множество~$\Gamma$ состоит из замыканий конечного числа критических траекторий
квадратичного дифференциала; поэтому множество регулярных точек компакта $\Gamma$ совпадает с
множеством $\mint{\Gamma}$ -- совокупностью открытых аналитических дуг,
составляющих~$\Gamma$.}.
Точнее,
\begin{lemma}\label{lemma2}
Пусть $F\in\kk_f$ и $(E,F,\psi)\in\sN_\psi$. Тогда для каждой
точки $z_0\in\mint F$ существует $*$-симметричная окрестность~$\sU\ni z_0$
такая, что
\begin{equation}
(3V^\lambda+G^\lambda_E+\psi)(z)=
(3V^\lambda+G^\lambda_E+\psi)(z^{*}),\qquad z\in\sU.
\label{s1.0}
\end{equation}
\end{lemma}

\begin{proof}
Действительно, пусть
$\pfi\colon \sU\to\DD$ -- отображение, фигурирующее в определении операции~$*$ (для
данной точки~$z_0\in\mint\Gamma$). Положим
\begin{equation}
\begin{gathered}
g(\zeta)=(3V^\lambda+G^\lambda_E+\psi-\gamma)(\pfi^{-1}(\zeta)),\qquad \zeta\in\DD,
\\
u(\zeta)=
\begin{cases}
\hphantom{-}g(\zeta),& \Im\zeta\geq0,\\
 -g(\zeta),& \Im\zeta<0,
\end{cases}
\end{gathered}
\label{ss1.01}
\end{equation}
где $\gamma=\gamma(\Gamma,\psi)$. Поскольку
$3V^\lambda+G^\lambda_E+\psi-\gamma\equiv0$ на $F\cap\sU$,
функция $u$ равна нулю на $(-1,1)$ и непрерывна в круге~$\DD$. Из
$\sS$-свойства~\eqref{a4} вытекает, что $u$ -- гладкая функция в~$\DD$; тем
самым, $u$ -- гармоническая функция в~$\DD$. Из принципа симметрии для
гармонических функций имеем $u(\zeta)=-u(\myo\zeta)$; возвращаясь в~$\sU$,
получаем~\eqref{ss1.01}.
\end{proof}

\begin{propos}[(аналог леммы 9 из~\cite{GoRa87})]\label{pro4}
Пусть мера $\mu$, $|\mu|\le1$, симметрична относительно вещественной оси,
отлична от равновесной меры $\lambda=\lambda(E;F,\psi)$
задачи~\eqref{a3}--\eqref{a4} и
пусть компакт $e\subset F$, $\mcap e=0$. Тогда существует точка $z_0\in F^0$, ${*}$-симметричная
окрестность $\sU$ точки $z_0$ и мера $\sigma$ такие, что:
\begin{enumerate}
\item[\rm (i)]
$|\sigma|<1$, $S(\sigma)\subset F\cup\sU$,
\item[\rm (ii)]
$\sigma\bigr|_{\sU}=\(\mu\bigr|_{\sU}\)^{*}$,
\end{enumerate}
и для функции
\begin{align}
V(z):
&=2V^\mu(z)+V^\sigma(z)+G_E^\sigma(z)+\psi(z) \notag\\
&=2V^\mu(z)+2V^\sigma(z)-V^{\myt\sigma}(z)+\psi(z)+C,
\label{4.2}
\end{align}
где $\myt\sigma$ -- выметание меры $\sigma$ на пластину $E$ конденсатора
Наттолла, имеем:
\begin{enumerate}
\item[\rm (iii)]
$V(z)=+\infty$, $z\in e$;
\item[\rm (iv)]
$V(z^*)=V(z)$, $z\in\sU$;
\item[\rm (v)]
$V(z_0)<V(z)$, $z\in F\setminus\{z_0\}$.
\end{enumerate}
\end{propos}
Утверждение~\ref{pro4} -- аналог леммы~9 из~\cite{GoRa87} (ср.~\cite[лемма
3]{RaSu13}).


\subsection{}\label{s4s4}
Итак, наши дальнейшие рассуждения основаны на следующем соотношении
(см.~\eqref{4.15.2}):
\begin{equation}
0=
\int_{\myFF}Q_{n,2}(t)
\biggl\{q_{m}(t)\myh{\rho}_n(t)-p_m(t)\biggr\}R_{n+m}(t)f^{*}(t)
\,dt,
\label{5.1}
\end{equation}
где $p_m,q_m$ -- {\it произвольные}
полиномы с комплексными коэффициентами, $p_m,q_m\in\PP_{m-1}$.

{\bad
Воспользуемся теперь следующей (предварительной) асимптотической формулой
($\LG$-приближением),
полученной Дж.~Наттоллом~\cite[формулы~(4.6)--(4.13)]{Nut86} для функции остатка
с помощью дифференциального уравнения Лагерра~\eqref{b1.30}
(см.~\cite[уравнение~(2.4)]{Nut86}, а также~\cite{MaRaSu12},~\cite{Rak12}),
которому удовлетворяет функция остатка и полиномы $P_{n,0}$ (см.~\eqref{01.7}):
}
\begin{equation}
R_n(z)=\biggl\{\frac{f(z)(\myt z_n-z)}{A_3(z)}\biggr\}^{1/2}
F_n(z)^{-1/4}\exp\{-n\psi_n(z)\}\{1+h_n(z)\},
\label{5.11}
\end{equation}
где
\begin{equation}
\psi_n(z)=\int_{z^0}^z\biggl\{\frac{(t-v_n)(t-\myt z_n)}{A_3(t)(t-z_n)}
\biggr\}^{1/2}\,dt.
\label{5.12}
\end{equation}
{\bad
Асимптотическая формула~\eqref{5.11} получена в~\cite{Nut86} классическим
методом Лиувилля--Грина и представляет собой
так называемое $\LG$-приближение
(см.~\cite{Olv74},~\cite[глава 2, \S~2, п.~2.1]{Fed86})
для решения уравнения Лагерра. В~\cite{Nut86} доказано, что
}
\begin{equation}
v_n=v+O(n^{-2/3}),\quad |\myt z_n-z_n|\leq\const(|z_n|+1)n^{-2/3},
\qquad n\to\infty,
\label{5.13}
\end{equation}
и приведен явный вид функции $F_n(z)$ (см.~\cite[(3.10)--(3.15)]{Nut86}).
В~\eqref{5.11} интегрирование проводится по ``поступательному пути'',
а динамика точи $\myt z_n$ при $n\to\infty$ подчиняется задаче обращения Якоби
(точные определения см.~\cite{Nut86}),

Пусть
\begin{equation}
\myh{\rho}_n(z):=
\int\frac{d\rho_n(\zeta)}{z-\zeta},
\label{5.14}
\end{equation}
тогда соотношение~\eqref{4.15.2} перепишется в виде
\begin{equation}
0=
\int_{\myFF}Q_{n,2}(t)
\bigl\{q_{m}(t)\myh{\rho}_n(t)-p_m(t)\bigr\}
R_{n+m}(t)f^{*}(t)\,dt,
\label{5.15}
\end{equation}
где $q_m,p_m$ -- {\it произвольные полиномы}, $q_m,p_m\in\PP_{m-m_0}$, $m=n/2$,
$m_0\geq2$ -- некоторое фиксированное число.

Зафиксируем теперь произвольный компакт $K\subset D(S)=\myo{\CC}\setminus{S}$
и произвольную последовательность монических полиномов $\omega_n(z)$ степени
$\leq{n-2m_0}$ с нулями на~$K$. Полиномы $q_m,p_m\in\PP_{m-m_0}$ в
соотношении~\eqref{5.15} произвольны, поэтому их можно выбрать так, чтобы
$\mdeg{q_m},\mdeg{p_m}\leq\mdeg{\omega_n}/2-1$ и
выполнялось условие интерполяции:
\begin{equation}
\frac{(q_m\myh{\rho}_n-p_m)(z)}{\omega_n(z)}\quad\text{-- голоморфная функция
на~$K$, $n=1,2,\dots$.}
\label{5.16}
\end{equation}
Тогда в силу~\eqref{5.14} и~\eqref{5.16} полиномы $q_m(z)=q_m(z;\omega_n)$
будут ортогональны по мере $d\rho_n$ с весом $1/\omega_n$, т.е. имеют место
соотношения
\begin{equation}
\int_E q_m(\zeta)\zeta^k\frac{d\rho_n(\zeta)}{\omega_n(\zeta)}=0,
\qquad k=0,1,\mdeg{\omega_n}/2-1,
\label{5.17}
\end{equation}
а соотношения~\eqref{5.15} эквивалентным образом переписываются в следующем
виде:
\begin{equation}
0=\int_{\myFF}Q_{n,2}(t)\frac{\omega_n(t)}{q_m(t)}
\biggl\{\int_{\EE}\frac{q^2_m(u)\,d\rho_n(u)}{\omega_n(u)(t-u)}
\biggr\}
R_{n+m}(t)\,f^{*}(t)\,dt,
\label{5.18}
\end{equation}
где теперь уже $\omega_n$ -- {\it произвольный полином} степени $\leq{n-2m_0}$
с нулями на компакте~$K$, $K$ -- произвольный зафиксированный компакт,
$K\subset D(S)$, $K\cap{\EE}=\varnothing$. При этом из леммы~\ref{le1} вытекает,
что
$$
d\rho_n(\zeta)\to f_j(\zeta)\,d\zeta,
\quad n\to\infty,\quad \zeta\in\EE^0_j,
\quad j=1,2,3.
$$
\end{proof}

Теперь, если
$$
\mu_n:=\frac1n\mu(Q_{n,2})\not\to\lambda_{\myFF},
$$
по некоторой подпоследовательности $n\in\Lambda$ имеем
\begin{equation}
\mu_n:=\frac1n\mu(Q_{n,2})\to\mu\neq\lambda_{\myFF},
\label{5.19}
\end{equation}
$|\mu|\leq1$. Опираясь на пару соотношений~\eqref{5.17} и~\eqref{5.18},
применим теперь с помощью утверждения~\ref{pro4} общий метод
Гончара--Рахманова, позволяющий находить асимптотику интегралов
вида~\eqref{6.2}, а затем и интеграла такого вида, как стоящего в правой
части~\eqref{5.18}. Предположим, как обычно~\cite{GoRa87},
что по некоторой подпоследовательности $n\in\Lambda$ имеем
\begin{equation}
\mu_n:=\frac1n\mu(Q_{n,2})\to\mu\neq\lambda_{\myFF}.
\label{5.19}
\end{equation}
Тогда в соответствии с утверждением~\ref{pro4} найдется мера~$\sigma$ такая, что
$|\sigma|<1$ и имеет место аналог леммы~9 из~\cite{GoRa87}. Выберем полином
$\omega_n$ так, чтобы
\begin{equation}
\frac1n\mu(\omega_n)\to\sigma,\quad n\to\infty,
\label{5.20}
\end{equation}
и воспользуемся соотношением ортогональности~\eqref{5.18} для того, чтобы вполне аналогично
работе~\cite{GoRa87} придти к противоречию. Полученное противоречие доказывает,
что
$$
\mu_n:=\frac1n\mu(Q_{n,2})\to\lambda_{\myFF},\quad n\to\infty.
$$


\section{Случай полиномов Паде}\label{s5}

\subsection{}\label{s5s1}
Изложим кратко основные положения нашей эвристической теории для изучения
глобальной структуры линий Стокса в применении
к полиномам Эрмита--Паде для пары функций $1,f$.
Обсудим вопрос о {\it сильной} асимптотике полиномов Паде
(см.~\eqref{5.11}--\eqref{5.12}).
Хорошо  известно (см.~\cite{Lag85},~\cite{Chu80},~\cite{Nut86}, а
также~\cite{Rak12},~\cite{MaRaSu13}),
что функция остатка $T_n$, произведение $P_{n,1}f$ и полином Паде $P_{n,0}$
удовлетворяют следующему дифференциальному уравнению Лагерра:
\begin{equation}
A_3(z)\Pi_1(z)w''+\Pi_3(z)w'+\Pi_2(z)w=0,
\label{b1.30}
\end{equation}
где $A_3(z)=\prod_{j=1}^3(z-a_j)$, $\Pi_k\in\CC_k[z]$, $k=1,2,3$, -- полиномы с
комплексными коэффициентами степени ровно\footnote{Мы считаем, что находимся
в ``общем положении''. Это в частности означает, что там, где надо считать,
что для некоторого полинома $p$ фиксированной степени $\mdeg{p}\leq{m}$, мы
считаем, что $\mdeg{p}={m}$.}~$k$. Точнее, имеем
\begin{equation}
\begin{gathered}
\Pi_1(z)=z-z_n,\qquad \Pi_2(z)=-n(n+1)(z-b_n)(z-v_n),\\
\Pi_3(z)=(z-z_n)B_2(z)-A_3(z),\qquad B_2(z)=A'_3(z)-A_3(z)\frac{f'}{f}(z).
\end{gathered}
\label{b1.31}
\end{equation}
Таким образом,~\eqref{b1.30} -- это дифференциальное уравнение с
полиномиальными коэффициентами и большим параметром при свободном члене.
Используя теоремы Шталя, Наттолл~\cite{Nut86} показал, что линиями Стокса
для этого дифференциального уравнения~\eqref{b1.30} являются в точности три
дуги $S_1,S_2,S_3$ компакта Чеботырёва--Шталя -- критические траектории
квадратичного дифференциала
$$
-\frac{z-v}{A_3(z)}\,dz^2.
$$
Сопряженные траектории этого квадратичного дифференциала являются
каноническими путями для соответствующей $\LG$-теории. Это наблюдение позволило
Наттоллу найти $\LG$-приближения для уравнения~\eqref{b1.30} и в итоге получить
формулы сильной асимптотики для функции остатка $R_n$ и полиномов Паде
$P_{n,j}$, $j=0,1$.
А именно, в~\cite{Nut86} доказано, что в~\eqref{b1.31} $v_n\to{v}$,
$\dist(z_n,b_n)\to0$,
$n\to\infty$, точка $z_n$ является решением специальной (эллиптической)
проблемы обращения Якоби.

В настоящей работе мы выведем формулы сильной асимптотики для функции остатка
и полиномов Паде с помощью нашего эвристического подхода, эти формулы следует
понимать в смысле сходимости по емкости, они будут не такие детальные как
формулы Наттолла~\eqref{5.11}. Однако сам вывод будет короче, чем вывод, данный
в его работе~\cite{Nut86}, поскольку используется тот факт, что соотношения
о сходимости по емкости~\eqref{01.12}--\eqref{01.14} можно дифференцировать
(см. лемму~\ref{lem2}).

Начнем свои рассуждения с вывода дифференциального уравнения~\eqref{b1.30}
(подробнее см.~\cite{Nut86},~\cite{MaRaSu12}).
Для того, что не усложнять обозначения, будет здесь опускать индекс $n$.
Поскольку функция $f$ имеет вид~\eqref{01.1}, то функция остатка $w=T_n(z)$ --
многозначная аналитическая функция в $\CC':=\myo\CC\setminus\{a_1,a_2,a_3\}$.
Нетрудно увидеть, что можно выбрать две независимые ветви этой функции $w_1=T_1$
и $w_2=T_2$ (или взять $P_{n,0}$ и $P_{n,1}f$ такие, что любая другая ветвь $R$ может быть записана в виде их
линейной комбинации:
$$
T(z)=C_1T_1(z)+C_2T_2(z),\qquad z\in\CC';
$$
например, можно взять $T_1=P_{0}$, $T_2=fP_{1}$. Тогда имеем:
$$
T'(z)=C_1T'_1(z)+C_2T'_2(z),\quad T''(z)=C_1T''_1(z)+C_2T''_2(z).
$$
Следовательно,
\begin{equation}
\det\begin{pmatrix}
T&T'&T''\\
T_1&T'_1&T''_1\\
T_2&T'_2&T''_2
\end{pmatrix}\equiv0.
\label{b1.32}
\end{equation}
При этом, с одной стороны, при $T_1=P_{0}$, $T_2=fP_{1}$ коэффициенты
уравнения~\eqref{b1.32} -- полиномы от $z$. С другой стороны, если к третьей
строке прибавить вторую, то в третьей строке получим слагаемые порядка $O(z^{-n-1}),
O(z^{-n-2}),O(z^{-n-3})$ при $z\to\infty$. Непосредственно отсюда и приходим
к уравнению~\eqref{b1.30}.

\subsection{}\label{s5s2}
К асимптотическим формулам типа~\eqref{5.11} можно придти и непосредственно из
уравнения~\eqref{b1.32}, если (после прибавления второй строки к третьей),
воспользоваться вытекающими из формул~\eqref{01.12}--\eqref{01.14}
предельными соотношениями (см. лемму~\ref{lem2})
\begin{gather}
\frac1n\frac{P'_{n,j}(z)}{P_{n,j}(z)}\overset\mcap\longrightarrow
-\sV'_\lambda(z)=\myh\lambda(z),
\quad n\to\infty,\quad z\in D,\quad j=0,1,\notag
\\
\frac1n\frac{R'_{n}(z)}{R_{n}(z)}\overset\mcap\longrightarrow
\pfi'(z)=-\myh\lambda(z),
\quad n\to\infty,\quad z\in D;
\notag
\end{gather}
здесь $\pfi(z)=-P(z,\infty)$, $P(z,\infty)=g(z,\infty)+ig^{*}(z,\infty)$ --
комплексная функция Грина. Напомним, что $g(z,\infty)=\gamma_S-V^\lambda(z)$,
\begin{equation}
\pfi'(z)=-G'(z,\infty)=-\sqrt{\frac{z-v}{A_3(z)}}.
\label{b1.36}
\end{equation}
Используя~лемму~\ref{lem2} непосредственно из~\eqref{b1.32}
получаем следующее предельное представление для соответствующего
характеристического уравнения:
\begin{equation}
\det\begin{pmatrix}
1&p&p^2\\
1&\pfi'&(\pfi')^2\\
1&-\pfi'&(\pfi')^2
\end{pmatrix}\equiv0.
\label{b1.37}
\end{equation}
Два независимых решения~\eqref{b1.37} это $p_{1,2}(z)=\pm\pfi'(z)$, при этом
из~\eqref{b1.36} вытекает, что $\pfi'(z)\neq-\pfi'(z)$, $z\notin{S}$.
Тем самым дифференциальное уравнение~\eqref{b1.30} принимает вид
\begin{equation}
w''+\delta_{n,1}(z)nw'-n^2(\pfi'(z))^2(1+\delta_{n,2}(z))w=0,
\label{b1.38}
\end{equation}
где $\delta_{n,j}\overset\mcap\longrightarrow0$, $n\to\infty$, $j=1,2$
(сходимость -- по емкости и с геометрической скоростью). Масштабирующий
множитель $n^2$ при свободном члене $w$ возник в уравнении~\eqref{b1.38} из-за
коэффициентов $1/n^2$ в соотношениях~\eqref{01.16} и~\eqref{01.17};
множителем $n$ при $w'$ можно пренебречь, поскольку $\delta_{n,1}\to0$
геометрически.

Тем самым в качестве
двух $\LG$-приближений $\myt{w}_{n,j}$ естественно рассматривать два решения
дифференциального уравнения
\begin{equation}
\myt w''-n^2(\pfi'(z))^2\myt w=0,
\label{b1.39}
\end{equation}
эти два решения имеют вид
\begin{equation}
\myt{w}_{n,j}(z)
=\frac1{\sqrt{n\pfi'(z)}}\exp\biggl\{\pm n\int_{a_1}^z\pfi'(\zeta)
\,d\zeta\biggr\}
=\frac1{\sqrt{n\pfi'(z)}}e^{\pm n\pfi(z)},\qquad
j=1,2.
\label{b1.40}
\end{equation}
Из~\eqref{b1.40} получаем следующие асимптотические формулы для нормированных
функции остатка и полиномов Паде
\begin{equation}
\begin{gathered}
(-1)^jP_{n,j}f^{j}\overset\mcap=
\frac1{\sqrt{\pfi'(z)}}e^{-n\pfi(z)}(1+o(1)),\qquad j=0,1,\\
R_{n}\overset\mcap=\frac1{\sqrt{\pfi'(z)}}e^{n\pfi(z)}(1+o(1)).
\end{gathered}
\label{b1.41}
\end{equation}
Из~\eqref{b1.41} вытекает, что
$$
(f-[n/n]_f)(z)\overset\mcap=
e^{2n\pfi(z)}(1+o(1)).
$$

\section{Построение $\LG$-приближений с помощью трехлистной
римановой поверхности с каноническим разбиением на листы}\label{s6}

{\bad
Опишем теперь соответствующую эвристическую процедуру построения
$\LG$-приближения для полиномов Эрмита--Паде. Для этого сначала методом
Гончара (см.~\cite[п.~3]{GoRaSu91}) с помощью симметричного конденсатора
$(E,F)$ построим каноническую трехлистную риманову поверхность.
}

\subsection{Построение трехлистной римановой поверхности с каноническим
разбиением на листы}\label{s6s1}

Построение трехлистной римановой поверхности $\RS_3$ основано на существовании
конденсатора Наттолла $(E,F)$ (см. утверждение~\ref{pro1}). Принцип построения
(т.е. нумерация листов и взаимные переходы с одного листа на
другой) такой трехлистной римановой поверхности был предложен А.~А.~Гончаром
в 1991 году~\cite[п.~3]{GoRaSu91} (в этой совместной работе именно Гончару принадлежат
формулы~\eqref{ago1}, непосредственно из которых вытекает наличие на~$\RS_3$
{\it комплексной структуры}; нужные нам здесь
формулы~\eqref{ago2} являются непосредственным следствием~\eqref{ago1}; см.
также~\cite{GoRaSu92} и~\cite[п.~6]{MaRaSu11b}). Отметим, что этот
принцип оказывается справедливым и в случае несвязного множества
$D^{*}=\myo\CC\setminus{F}$ (см. в~\cite{RaSu13} некоторые частные примеры).

Приведем явное в терминах конденсатора Наттолла $(E,{F})$ и равновесных
мер $\lambda=\lambda_E$ и $\myt\lambda=\lambda_{F}$ (напомним, что
$\myt{\lambda}=\beta_F(\lambda_E)$ -- выметание меры $\lambda_E$ из $D^*$ на~$F$)
описание канонической трехлистной римановой поверхности.

Построим трехлистную риманову поверхность
$\RS_3=\myo{\RS^{(1)}}\cup\myo{\RS^{(2)}}\cup\myo{\RS^{(3)}}$, состоящую из
трех замкнутых листов, следующим образом. Возьмем три
экземпляра римановой сферы $\myo{\CC}$. На первом проведем разрез по дугам
$E_j$, $j=1,2,3$, на втором
-- по $E_j$ и дугам $F_k$, $k=1,\dots,6$, составляющим ${F}$, на третьем --
только по дугам $F_k$ (см. рис.~\ref{Fig3}).
Три полученных экземпляра римановой сферы переклеиваются друг с
другом следующим образом. Второй подклеивается к первому по разрезу,
соответствующему дугам~$E_j$, третий -- ко второму по разрезам,
соответствующим дугам $F_k$. Полученная трехлистная риманова поверхность
$\RS_3=\myo{\RS^{(1)}}\cup\myo{\RS^{(2)}}\cup\myo{\RS^{(3)}}$
гомеоморфна римановой сфере с одной ручкой (род $\RS_3$ равен $g=1$).
Определим на $\RS_3$ функцию $u(\zz)$, $\zz\in\RS_3$, следующим образом:
\begin{equation}
u(z^{(1)})=2G^{\lambda}_{{F}}(z),\quad
u(z^{(2)})=G^{\lambda}_{{F}}(z)-3V^{\lambda}(z)+\gamma,\quad
u(z^{(3)})=-G^{\lambda}_{{F}}(z)-3V^{\lambda}(z)+\gamma.
\label{ago1}
\end{equation}
Непосредственно из условия равновесия~\eqref{2.6} и $\sS$-свойства~\eqref{a4}
компакта ${F}$ вытекает, что $u$ -- гармоническая функция на
$\RS_3\setminus({F}^{(1)}\cup\EE^{(3)})\setminus\{\infty^{(2)},\infty^{(3)}\}$.
Кроме того, $u\equiv0$ на компакте $\KK=F^{(1)}\cup\EE^{(3)}$ и
$u(z^{(2)})=3\log|z|+c_2+o(1)$ при $z^{(2)}\to\infty^{(2)}$,
$u(z^{(3)})=3\log|z|+c_3+o(1)$ при $z^{(3)}\to\infty^{(3)}$.
Отсюда вытекает, что
$$
u(\zz)\equiv3(g_{\KK}(\zz,\infty^{(2)})+g_{\KK}(\zz,\infty^{(3)})),
$$
где
$g_{\KK}(\zz,\,\cdot\,)$ -- функция Грина для области $\RS_3\setminus{\KK}$ с
особенностью в соответствующей точке.
Теперь уже нетрудно проверить, что на построенной римановой поверхности $\RS_3$
следующая многозначная аналитическая функция является $\phi$-функцией Наттолла
(см.~\cite[\S~3, п.~3.1, формулы~(3.1.2)--(3.1.5)]{Nut84}:
\begin{equation}
\phi(z^{(1)})=2\sV^{\lambda}(z)-2\gamma,\quad
\phi(z^{(2)})=\sV^{\myt\lambda}(z)-2\sV^{\lambda}(z)-\gamma+c,\quad
\phi(z^{(3)})=-\sV^{\myt\lambda}(z)-\gamma+2c,
\label{ago2}
\end{equation}
где $c=\const$ -- некоторая постоянная (благодаря которой функция $\Re\phi$
склеивается через разрезы в непрерывную функцию на $\RS_3$),
$$
\sV^{\mu}(z)=-\int\log(z-\zeta)\,d\mu(\zeta)
$$
-- комплексный потенциал меры~$\mu$; $\phi(\zz)$ -- многозначный абелев
интеграл третьего рода на $\RS_3$. Нетрудно увидеть, что $\phi$-функция Наттолла,
введенная им еще в 1981 году в~\cite{Nut81}, совпадает с $g$-функцией
Дейфта~\cite{Dei99}.
Функция $\Re\phi(\zz)$ уже однозначная
гармоническая функция на
$\RS_3\setminus\{\infty^{(1)},\infty^{(2)},\infty^{(3)}\}$.
Из~\eqref{ago2} вытекает, что при $z\in\myo\CC\setminus(E\cup{F})$
справедливы неравенства
$$
\Re(\phi(z^{(3)})-\phi(z^{(2)}))>0,\qquad
\Re(\phi(z^{(2)})-\phi(z^{(1)}))>0.
$$
Тем самым, в полном соответствии с гипотезой
Наттолла~\cite[формула~(3)]{Nut82}, замкнутые листы римановой
поверхности $\RS_3$ оказываются упорядоченными следующим образом:
\begin{equation}
\Re\phi(z^{(3)})\geq\Re\phi(z^{(2)})\geq\Re\phi(z^{(1)}),\qquad z\in\myo\CC,
\label{anut1}
\end{equation}
причем равенство выполняется только на линиях склейки ${F}^{(3,2)}$
и~$E^{(2,1)}$, для открытых листов $\RS^{(j)}$, $j=1,2,3$, в~\eqref{anut1}
имеет место строгое неравенство, а $\RS_3$ представляется в виде
$\RS_3=\RS^{(3)}\sqcup{F}^{(3,2)} \sqcup\RS^{(2)}\sqcup
E^{(2,1)}\sqcup\RS^{(1)}$.
Функция $\Phi(\zz):=e^{\phi(\zz)}$ -- многозначная аналитическая функция на
$\RS_3$ с {\it однозначным модулем} и
$E=\{z\in\myo\CC:|\Phi(z^{(1)})|=|\Phi(z^{(2)})|\}$,
${F}=\{z\in\myo\CC:|\Phi(z^{(3)})|=|\Phi(z^{(2)})|\}$.
При этом $d\phi(\zz)$ -- абелев дифференциал третьего рода на $\RS_3$ с чисто
мнимыми периодами и особенностями в точках
$\zz=\infty^{(1)},\infty^{(2)},\infty^{(3)}$ вида $2/z,-1/z,-1/z$ соответственно
(сумма всех вычетов равна нулю). Равновесные меры $\lambda,\myt\lambda$ --
единичные меры с носителями на $E=\bigcup_{j=1}^3E_j$ и
${F}=\bigcup_{k=1}^6{F}_k$ соответственно.

Функция
$$
y(\zz)=\phi'(\zz)=\frac{d\phi}{dz}(\zz)
$$
будет уже {\it однозначной} мероморфной функцией на $\RS_3$.
Положим $y_j(z):=y(z^{(j)})$,
$j=1,2,3$.
Тогда $y_j$ -- многозначные функции относительно
$z\in\myo\CC'':=\myo\CC\setminus\{a_1,a_2,a_3,v_1,v_2\}$.
Имеем
\begin{equation}
\det
\begin{pmatrix}
y^3&y^3_1&y^3_2&y^3_3\\
y^2&y^2_1&y^2_2&y^2_3\\
y&y_1&y_2&y_3\\
1&1&1&1
\end{pmatrix}
\equiv0,\qquad z\in\myo\CC''.
\label{a5.3}
\end{equation}
Из~\eqref{a5.3} вытекает, что функция $y$ удовлетворяет
следующему алгебраическому уравнению 3-го порядка (кубическому уравнению,
ср.~\cite{NuTr87},~\cite[формулы~(2.1)--(2.8)]{ApKu11}):
\begin{equation}
y^3+r_2(z)y^2+r_1(z)y+r_0(z)\equiv0,\quad z\in\myo\CC',\quad r_j\in\CC(z),
\label{a5.4}
\end{equation}
$r_j\in\CC(z)$, $j=1,2,3$, -- рациональные функции от~$z$. Действительно,
пусть $A_j(z)$, $j=1,2,3,4$, -- алгебраические дополнения к элементам первого
столбца в~\eqref{a5.3}, тогда
$$
\sum_{j=1}^4 A_j(z) y^{j-1}\equiv0.
$$
Нетрудно увидеть, что функции $A_j(z)/A_4(z)$, $j=1,2,3$, -- однозначные
голоморфные функции переменного $z$ с возможными особенностями типа полюса
только в точках $\{a_1,a_2,a_3,v_1,v_2\}$. Следовательно,
$A_j(z)/A_4(z)\in\CC(z)$, $j=1,2,3$, -- рациональные функции и мы приходим
к~\eqref{a5.4}.

Отметим, что для открытых листов $\RS^{(j)}$, $j=1,2,3$, справедливы
соотношения:
$\proj{\RS^{(3)}}=D^{*}=\Omega({F})$,
$\proj{\RS^{(2)}}=\Omega({F})\setminus{E}$,
$\proj{\RS^{(3)}}=\Omega(E)=\myo\CC\setminus{E}$.
Совокупность кривых $E\cup{F}$ состоит из замыканий критических траекторий
квадратичного дифференциала
$$
-(\phi'(z)\,dz)^2>0.
$$


\subsection{Вывод дифференциального уравнения третьего порядка}\label{s6s2}

Начнем свои рассуждения с вывода дифференциального уравнения третьего порядка
(см.~\cite{MaRaSu13}).
Мы будем следовать здесь подходу Римана~~\cite{Rie68} (см. также~\cite{Chu80}
и~\cite{Nut86}).

Пусть $f$ -- голоморфная в бесконечно удаленной точки $z=\infty$ функция,
заданная представлением
\begin{equation}
f(z)=f(z;\balpha):=\prod_{j=1}^p(z-a_j)^{\alpha_j},\quad
\text{где}\quad
p\geq2,\quad\alpha_j\in\CC\setminus2\ZZ,\quad
\sum_{j=1}^p\alpha_j=0,
\label{ma1.1}
\end{equation}
и нормированная условием: $f(\infty)=1$. Функция $f$ продолжается как
многозначная аналитическая функция в область
$\Omega_{\ba}:=\myo\CC\setminus\{a_1,\dots,a_p\}$. В дальнейшем мы считаем, что
точки $a_1,\dots,a_p$ находятся в ``общем положении'' (см. ниже п.~\ref{s1s2});
при фиксированном $\ba=\{a_1,\dots,a_p\}$
класс функций вида~\eqref{ma1.1} будем обозначать через $\sL_{\ba}$,
$\balpha=\{\alpha_1,\dots,\alpha_p\}$.

Следующее дифференциальное уравнение второго порядка для полиномов Паде
степени $n$ функции $f$ (определение см. ниже п.~\ref{s1s2})
принадлежит Лагерру~\cite{Lag85} и хорошо известно
(см.~\cite[п.~5.1]{Nut84},~\cite{Nut86},~\cite{MaRaSu12},~\cite{Rak12}):
\begin{equation}
Ahw''+\{(A'-B)h-Ah'\}w'-n(n+1)gw\equiv0,
\label{ma1.5}
\end{equation}
где полиномы
\begin{equation}
A=A_p(z)=\prod\limits_{j=1}^p(z-a_j)\in\PP_p,
\qquad
B=B_p(z)=Af'/f\in\PP_{p-2}
\label{ma1.3}
\end{equation}
известны, а
$h=h_n(z)=z^{p-2}+\dotsb\in\PP_{p-2}$, $g=g_n=z^{2p-4}+\dotsb\in\PP_{2p-4}$ --
некоторые зависящие от $n$ полиномы
{\it фиксированной степени}~$\leq{2p-4}$ (здесь и всюду в
дальнейшем для произвольного $m\in\NN$ через $\PP_m$ обозначается множество
всех полиномов с комплексными коэффициентами степени $\leq{m}$; для
произвольного полинома $P\in\CC[z]$ через $P^{*}=P^{*}(z)=z^{\mdeg{P}}+\dotsb$
обозначается соответствующий монический полином).

Для частного случая, когда $p=2$ и
$f(z)=f(z;\alpha)=(z+1)^{-\alpha}(z-1)^{\alpha}=\bigl((z-1)/(z+1)\bigr)^{\alpha}$,
$\alpha\in\CC\setminus\ZZ$,
из~\eqref{ma1.5} вытекает хорошо известное (см.~\cite[п.~4.2]{Sze62})
дифференциальное уравнение для классических полиномов Якоби $P^{(\alpha,-\alpha)}_n(z)$:
\begin{equation}
(z^2-1)w''+2(z-\alpha)w'-n(n+1)w=0.
\label{ma1.3.1}
\end{equation}

Основной
результат работы~\cite{MaRaSu13}
-- следующее однородное дифференциальное
уравнение 3-го порядка, которому удовлетворяют полиномы Эрмита--Паде
$Q_{n,0},Q_{n,1},Q_{n,2}$ степени $n$ функции $f$:
\begin{equation}
A^2Hw'''+A\{3(A'-B)H-AH'\}w''-3(n-1)(n+2)Fw'+
2n(n^2-1)Gw\equiv0,
\label{ma1.2}
\end{equation}
где полиномы $A$ и $B$ имеют тот же смысл, что и выше,
$H=H_n(z)=z^{3p-6}+\dotsb\in\PP_{3p-6}$,
$F=F_n(z)=z^{5p-8}+\dotsb\in\PP_{5p-8}$ и
$G=G_n(z)=z^{5p-9}+\dotsb\in\PP_{5p-9}$ -- некоторые зависящие от $n$ полиномы
фиксированной степени.

Также как и для полиномов Паде, в частном случае, когда в~\eqref{ma1.1}
$p=2$ и $f(z)=f(z;\alpha)=\bigl((z-1)/(z+1)\bigr)^{\alpha}$,
$2\alpha\in\CC\setminus\ZZ$,
все полиномиальные коэффициенты степеней соответственно $4,3,2$
и $1$ уравнения~\eqref{ma1.2} найдены нами в {\it явном виде}, т.е. выражены
непосредственно в терминах номера $n$ и параметра $\alpha$ (см.~\eqref{ma1.9},
ср.~\eqref{ma1.3.1}).

{\bad
Уравнение~\eqref{ma1.2} естественно рассматривать как обобщение на случай
полиномов Эрмита--Паде дифференциального уравнения
Лагерра~\eqref{ma1.5}.
Для частного случая, когда $p=2$, соответствующие полиномы Эрмита--Паде
$Q^{(\alpha,-\alpha)}_{n,j}(z)$, $j=0,1,2$,
естественно рассматривать как обобщение классических ортогональных полиномов
Якоби~$P^{(\alpha,-\alpha)}_n(z)$.
}

Сформулируем основной результат работы~\cite{MaRaSu13}.
\begin{theorem}~\label{mat1}
Пусть функция $f(z)=f(z;\balpha)$ задана представлением~\eqref{ma1.1} и
нормирована условием $f(\infty)=1$, $Q_{n,j}$, $j=0,1,2$, -- полиномы
Эрмита--Паде степени $n$ и выполнены условия ``общего положения'' на точки
$\{a_1,\dots,a_p\}$. Тогда полиномы $Q_{n,0}$ и функции $Q_{n,j}f^{j}$,
$j=1,2$, являются независимыми
решениями следующего однородного дифференциального уравнения 3-й степени
с зависящими от $n$ полиномиальными коэффициентами фиксированной степени
$\leq{5p-6}$:
\begin{equation}
A^2Hw'''+A\{3(A'-B)H-AH'\}w''-3(n-1)(n+2)Fw'+2n(n^2-1)Gw\equiv0,
\label{ma1.8}
\end{equation}
где полиномы
$$
A=A_p(z)=\prod\limits_{j=1}^p(z-a_j)\in\PP_p,
\qquad
B=Af'/f\in\PP_{p-2},
$$
$H=H_n(z)=z^{3p-6}+\dotsb\in\PP_{3p-6}$,
$F=F_n(z)=z^{5p-8}+\dotsb\in\PP_{5p-8}$ и
$G=G_n(z)=z^{5p-9}+\dotsb\in\PP_{5p-9}$ -- некоторые зависящие от $n$
монические полиномы фиксированной степени.
\end{theorem}

\begin{remark}~\label{marem1}
Из определения~\eqref{1.5} вытекает, что функция остатка $R_n$ также
удовлетворяет уравнению~\eqref{ma1.8}.
\end{remark}

На основе теоремы~\ref{mat1} доказывается

\begin{theorem}~\label{mat2}
Пусть
$p=2$ и $f=f(z;\alpha)=\bigl((z-1)/(z+1)\bigr)^{\alpha}$,
где $2\alpha\in\CC\setminus\ZZ$. Тогда дифференциальное уравнение для полиномов
Эрмита--Паде $Q_{n,0}$ и функций $Q_{n,1}f$ и $Q_{n,2}f^2$ имеет следующий вид:
\begin{align}
(z^2-1)^2w'''
&+6(z^2-1)(z-\alpha)w''\notag\\
&-\bigl[3(n-1)(n+2)z^2+12\alpha z-(3n(n+1)+8\alpha^2-10)\bigr]w'\notag\\
&+2\bigl[n(n^2-1)z+\alpha(3n(n+1)-8)\bigr]w=0.
\label{ma1.9}
\end{align}
Если $\alpha\in(0,1/2)$, то все нули полиномов $Q_{n,j}$ лежат на
$\RR\setminus[-1,1]$ и их предельное распределение имеет плотность
\begin{equation}
\frac{\sqrt{3}}{2\pi}\,\frac{1}{\sqrt[3]{x^2-1}}\,\left(
\frac{1}{\sqrt[3]{|x|-1}}-\frac{1}{\sqrt[3]{|x|+1}}\right),\quad |x|>1.
\label{ma1.91}
\end{equation}
\end{theorem}


Анализировать это дифференциальное
уравнение~\eqref{ma1.2} столь же детально, как это было сделано Наттоллом
(см.~\cite{Nut86}, а также \S~\ref{s5}), пока не удается. Поэтому поступим
по другому, а именно также как в \S~\ref{s5}.

Используя формулы сходимости по емкости~\eqref{01.12}--\eqref{01.14}, их
сохранение при дифференцировании обеих частей (см. лемму~\ref{lem2}) и
найденное устройство трехлистной римановой
поверхности $\RS_3$, аналогично~\S~\ref{s5} получаем для
соответствующего~\eqref{ma1.2} характеристического уравнения
\begin{equation}
0\equiv\det\begin{pmatrix}
1&p&p^2&p^3\\
1&\pfi'(z^{(1)})&(\pfi')^2(z^{(1)})&(\pfi')^3(z^{(1)})\\
1&\pfi'(z^{(2)})&(\pfi')^2(z^{(2)})&(\pfi')^3(z^{(2)})\\
1&\pfi'(z^{(3)})&(\pfi')^2(z^{(3)})&(\pfi')^3(z^{(3)})
\end{pmatrix}
=\prod_{j=1}^3(p-\pfi'(z^{(j)})).
\label{b1.44}
\end{equation}
Таким образом, три независимых решения~\eqref{b1.44} это $p_j(z)=\pfi'(z^{(j)})$,
$j=1,2,3$, при этом $\pfi'(z^{(j)})\neq\pfi'(z^{(k)})$, $k\neq j$,
$z\notin{E\cup F}$.
Тем самым в качестве трех $\LG$-приближений $\myt{w}_{n,j}$ для решений
дифференциального уравнения~\eqref{ma1.2} естественно рассматривать три решения
дифференциального уравнения
\begin{equation}
\det\begin{pmatrix}
\myt{w}&\myt{w}'/n&\myt{w}''/n^2&\myt{w}'''/n^3\\
1&\pfi'(z^{(1)})&(\pfi')^2(z^{(1)})&(\pfi')^3(z^{(1)})\\
1&\pfi'(z^{(2)})&(\pfi')^2(z^{(2)})&(\pfi')^3(z^{(2)})\\
1&\pfi'(z^{(3)})&(\pfi')^2(z^{(3)})&(\pfi')^3(z^{(3)})
\end{pmatrix}
\equiv0
\label{b1.45}
\end{equation}
(сомножители $1/n,1/n^2,1/n^3$ появляются в~\eqref{b1.45} по той же причине,
что и в~\S~\ref{s5}).
Эти три решения имеют асимптотический вид, задаваемый формулой Федорюка
(см.~\cite[п.~3.2, п.~4.5 и формула~(2.50)]{Fed86})
\begin{equation}
\myt{w}_{n,j}(z)
\fallingdotseq
\exp\biggl\{n\pfi(z^{(j)})-\int_{a_1}^{z^{(j)}}\sum_{k=1,k\neq{j}}^3
\frac{\pfi''(z^{(j)})}{\pfi'(z^{(j)})-\pfi'(z^{(k)})}\,dz\biggr\},
\quad
j=1,2,3.
\label{b1.46}
\end{equation}
Так так риманова поверхность $\RS_3$ -- класса Наттолла, т.е.
\begin{equation}
\Re\pfi(z^{(3)})>\Re\pfi(z^{(2)})>\Re\pfi(z^{(1)}),\quad z\notin{E\cup F},
\label{b1.47}
\end{equation}
то $\LG$-приближения~\eqref{b1.46} асимптотически независимы при движении
по соответствующему каноническому пути (тем самым в открытом множестве
$\myo\CC\setminus(E\cup F)$ всегда выполняется
условие~(B) из~\cite[п.~3.1]{Fed86}).
Следовательно, аналогично классической формуле Биркгофа, любое решение
дифференциального уравнения~\eqref{ma1.2} может быть приближенно некоторой
линейной комбинацией
$$
C_1\myt{w}_{n,1}(z)+C_2\myt{w}_{n,2}(z)+C_3\myt{w}_{n,3}(z),
$$
в которой в силу неравенств~\eqref{b1.47} всегда есть только одно доминирующее
при $n\to\infty$ слагаемое. Выбор этого доминирующего слагаемого зависит
от того, какому именно листу рп $\RS_3$ соответствуют заданные начальные
условия для дифференциального уравнения~\eqref{b1.45}.

\section{Заключительные замечания}\label{s7}

\subsection{}\label{s7s1}
В работах М.~В.~Федорюка неоднократно отмечалось (см., например,
\cite[\S~3, п.~3.3; \S~4, п.4.5]{Fed86}), что до сих пор не существует
глобальной асимптотической теории для уравнений порядка $n\geq3$ и, по его
мнению, ее нельзя построить вообще.

Из сформулированных в настоящей работе результатов вытекает, что такую
глобальную асимптотическую теорию вполне можно построить в некотором классе
дифференциальных уравнений 3-го порядка, у которых соответствующая трехлистная
риманова поверхность допускает каноническое разбиение на
листы. Настоящая работа фактически
представляет собой {\it программу} такого построения с помощью метода
Гончара--Рахманова.

Схема применения метода Гончара--Рахманова для изучения распределения нулей
полиномов Эрмита--Паде и исследования глобальной топологии линий
Стокса для соответствующего дифференциального уравнения 3-го порядка состоит из
следующих шагов:

(1) правильная постановка экстремальной теоретико-потенциальной задачи
равновесия для семейства допустимых конденсаторов, находящихся во внешнем поле;

(2) доказательство того, что экстремальный  конденсатор $(E,F)$ обладает
$\sS$-свойством во внешнем поле, а его пластины $E$ и $F$ состоят из замыканий
критических траекторий квадратичного дифференциала;

(3) доказательство того, что распределение нулей полиномов Эрмита--Паде
совпадает с равновесной мерой для $F$-пластины экстремального конденсатора
и вывод формулы для асимптотики отношения этих полиномов;

(4) построение трехлистной римановой поверхности с каноническим разбиением
на листы;

(5) нахождения дифференциального уравнения 3-го порядка для полиномов
Эрмита--Паде и функции остатка; доказательство того, что на канонической
трехлистной римановой поверхности критические траектории квадратичного
дифференциала являются линиями Стокса для этого уравнения, а сопряженные
траектории -- каноническими путями для $\LG$-приближений решений этого
уравнения.

\subsection{}\label{s7s2}
Утверждается (в качестве гипотезы), что и в общем случае произвольной алгебраической функции $f$
построить конденсатор Наттолла для набора из трех функций $1,f,f^2$
невозможно без знания структуры компакта Шталя для функции $f$. А
именно, конденсатор Наттолла строится с помощью
функции Грина $g_{D}(z,\infty)$ для области Шталя $D$ и ``остова'' компакта
Шталя, который состоит из ``эффективных'' точек ветвления функции $f$ и
соответствующих им точек Чеботарёва. Точнее, пусть $-V(z)A(z)^{-1}\,dz^2$
-- квадратичный дифференциал, замыкания критических траекторий которого
составляют компакт Шталя~$S$, $A(z)=\prod_{j=1}^M(z-a_j)$, $\{a_,\dots,a_M\}$
-- точки ветвления функции~$f$. Положим $\myt{V}_{p-2}/\myt{A}_p\equiv
V_{p-2}/A_p$, где $\NOD(\myt{V}_{p-2},\myt{A}_p)=1$, т.е. полиномы
$\myt{V}_{p-2}$ и $\myt{A}_p$ взаимно
просты. Удалением всех сомножителей вида $(z-b_j)^{2s}$, $s\in\NN$, перейдем от
$\myt{V}_{p-2}$ к новому полиному, который будем также обозначать
$\myt{V}_{p-2}$. Тогда $\myt{V}_{p-2}\myt{A}_p$
-- полином четной степени без кратных корней.

\begin{definition}\label{eff_points}
Нули полинома $\myt{V}_{p-2}$ будем называть точками Чеботарёва для компакта
Шталя, а нули полинома $\myt{A}_p$ -- эффективными точками ветвления функции
$f$. Будем говорить, что множество точек Чеботарёва и эффективных точек
ветвления образуют остов компакта Шталя
$S_{\frame}=\{v_1,\dots,v_{N_1},a_1,\dots,a_{N_2}\}$.
\end{definition}

\subsection{Общий вывод}\label{s7s3}
Согласно сказанному выше, видимо имеет место следующий общий принцип,
справедливый, например, в классе алгебраических функций $f\in\HH(\infty)$:
\begin{enumerate}
\item[\rm (i)] для пары $1,f$ построим аппроксимации Эрмита--Паде (в данном
случае просто Паде) в точке $z=\infty$ и предельным переходом найдем компакт
Шталя, после этого построим остов Шталя $S_{\frame}$;
\item[\rm (ii)] для трех функций $1,f,f^2$
построим аппроксимации Эрмита--Паде в точке $z=\infty$ и предельным переходом
найдем конденсатор Наттолла;
\item[\rm (iii)] для набора функций $1,f,f^2,f^3$
построим аппроксимации Эрмита--Паде в точке $z=\infty$ и предельным переходом
найдем что-то пока неизвестное;
\item[\rm (iv)] $\dots$
\end{enumerate}

\begin{figure}[h!]
\centerline{
\includegraphics[width=13cm,height=17.5cm]{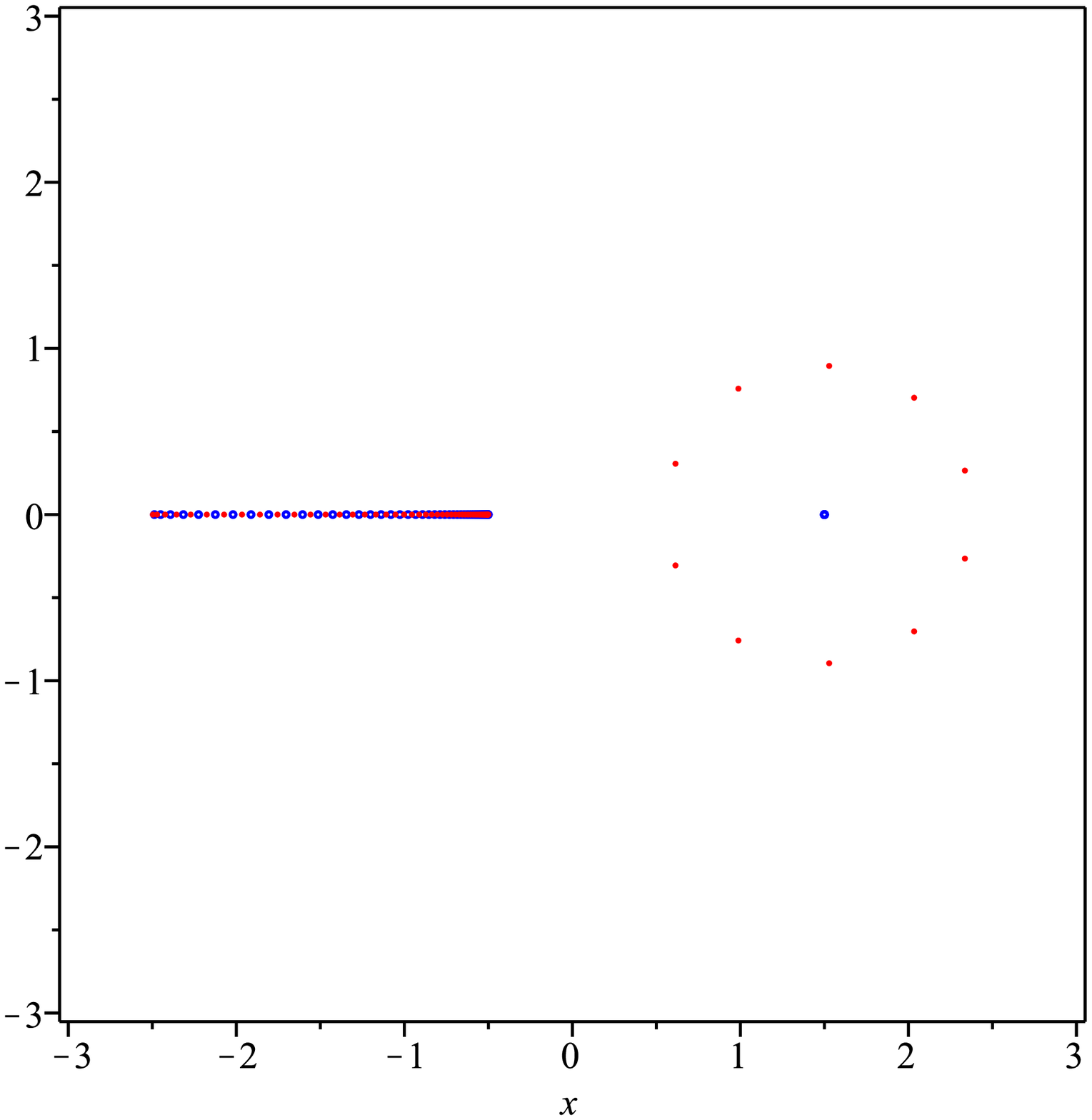}}
\vskip-6mm
\caption{Расположение линий Стокса на первом листе рп~$\RS_3$.}
\label{Fig1}
\end{figure}

\begin{figure}[h!]
\centerline{
\includegraphics[width=13cm,height=17.5cm]{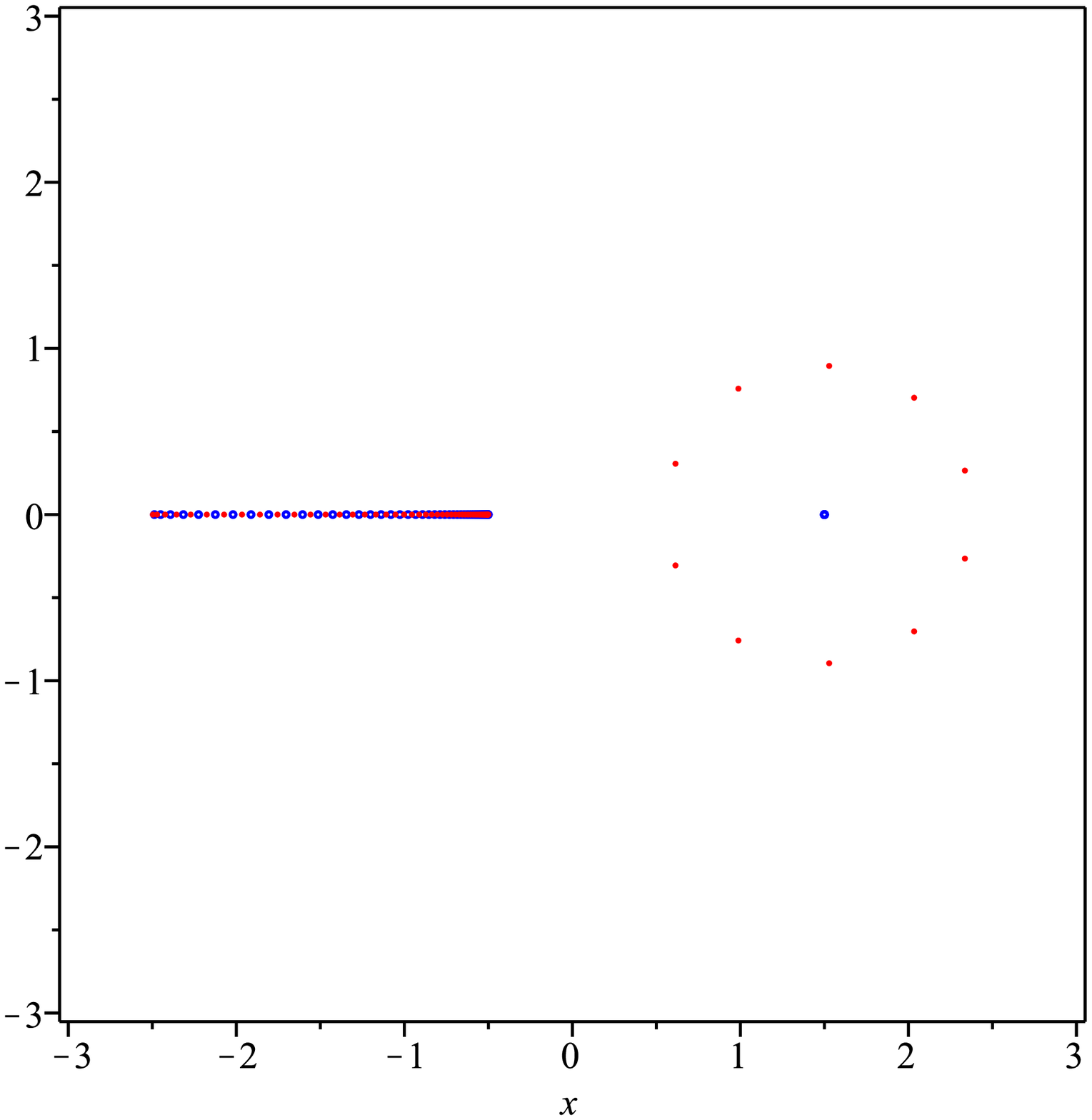}}
\vskip-6mm
\caption{Расположение линий Стокса на втором листе рп~$\RS_3$.}
\label{Fig2}
\end{figure}

\begin{figure}[h!]
\centerline{
\includegraphics[width=13cm,height=17.5cm]{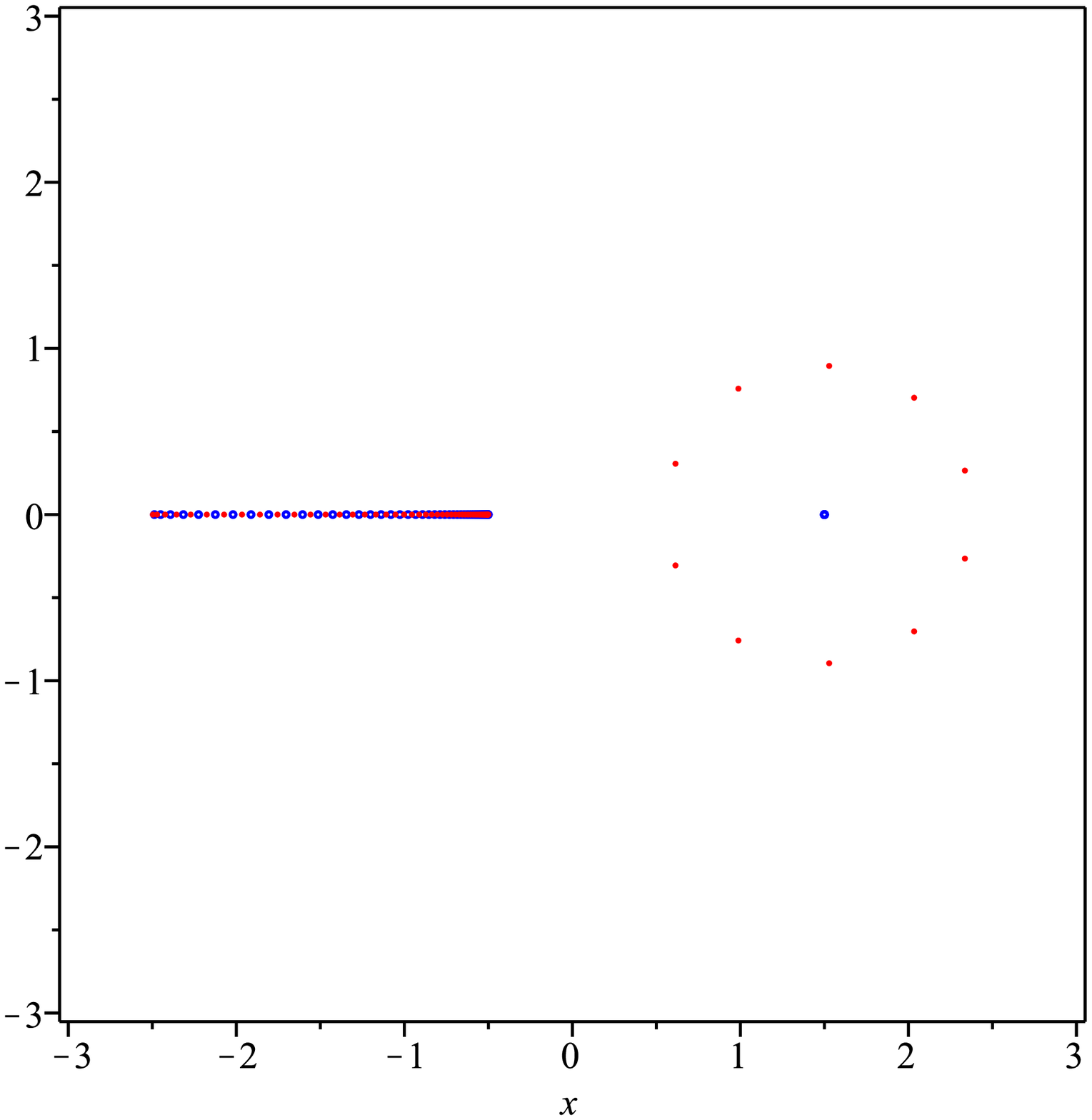}}
\vskip-6mm
\caption{Расположение линий Стокса на третьем листе рп~$\RS_3$.}
\label{Fig3}
\end{figure}

См. рис.~\ref{Fig1}--\ref{Fig3}.

\newpage

\newpage

\newpage

\newpage

\begin{figure}[h!]
\centerline{
\includegraphics[width=13cm,height=13cm]{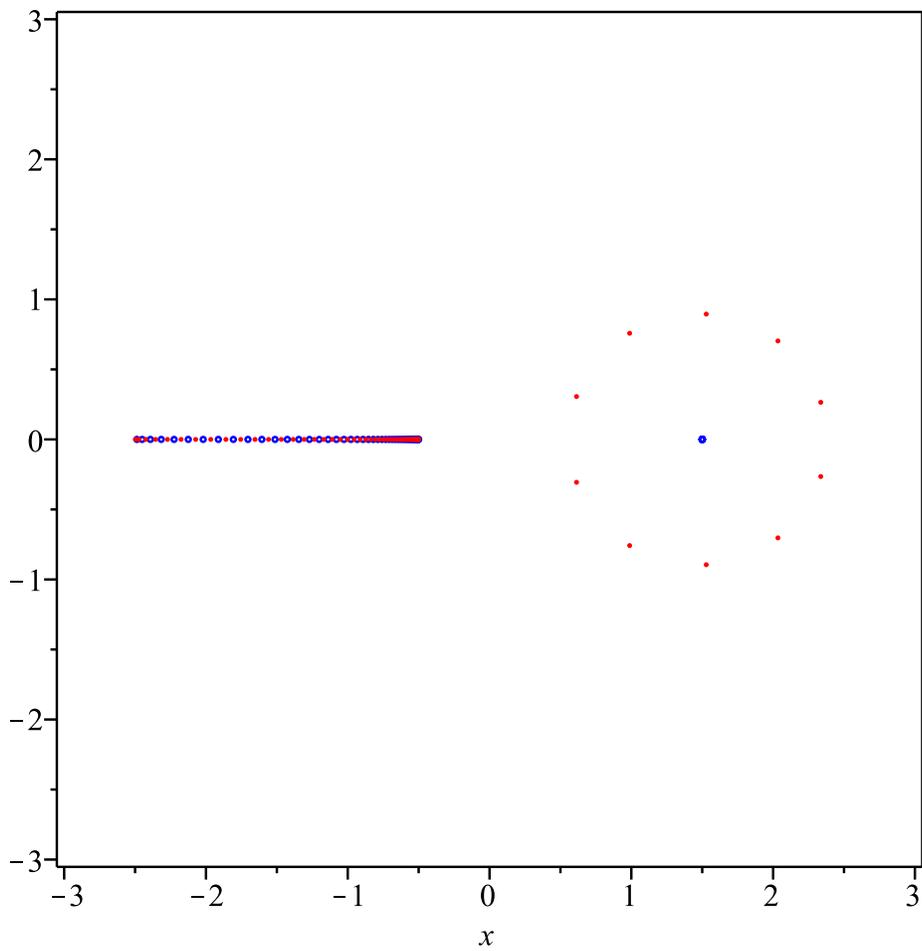}}
\vskip-6mm
\caption{Нули полиномов Паде для $[1,f,f^2]$,
$f=(1-z^2)^{1/3}(1-i\sqrt{3}z)^{-2/3}$.}
\label{Fig4}
\end{figure}
\newpage



\end{document}